\documentclass[12pt,oneside,reqno]{amsart}

\usepackage{amssymb}
\usepackage{txfonts}
\usepackage{bbm}
\usepackage{cases}
\usepackage{amsmath}
\usepackage{graphicx}
\usepackage{mathrsfs}
\usepackage{stmaryrd}
\usepackage{amsfonts}
\usepackage{enumerate,amsmath,amssymb,amsthm}

\numberwithin{equation}{section}

\newcommand{\be}{\begin{eqnarray}}
\newcommand{\ee}{\end{eqnarray}}
\newcommand{\ce}{\begin{eqnarray*}}
\newcommand{\de}{\end{eqnarray*}}
\newtheorem{theorem}{Theorem}[section]
\newtheorem{lemma}[theorem]{Lemma}
\newtheorem{remark}[theorem]{Remark}
\newtheorem{definition}[theorem]{Definition}
\newtheorem{proposition}[theorem]{Proposition}
\newtheorem{Examples}[theorem]{Example}
\newtheorem{corollary}[theorem]{Corollary}

\def\eps{\varepsilon}

\def\e{\mathrm{e}}

\def\v{\mathrm{v}}
\def\u{\mathbf{u}}
\def\w{\mathbf{w}}
\def\p{\partial}

\def\[{{\Big[}}
\def\]{{\Big]}}
\def\<{{\langle}}
\def\>{{\rangle}}
\def\({{\Big(}}
\def\){{\Big)}}

\def\bx{{\mathbf{x}}}
\def\tr{\mathrm {tr}}

\def\dif{{\mathord{{\rm d}}}}

\def\no{\nonumber}
\def\={&\!\!=\!\!&}

\def\cE{{\mathcal E}}

\def\cK{{\mathcal K}}

\def\cM{{\mathcal M}}

\def\cO{{\mathcal O}}
\def\cP{{\mathcal P}}

\def\cR{{\mathcal R}}

\def\mE{{\mathbb E}}

\def\mH{{\mathbb H}}
\def\mI{{\mathbb I}}

\def\mL{{\mathbb L}}

\def\mN{{\mathbb N}}

\def\mP{{\mathbb P}}
\def\mQ{{\mathbb Q}}
\def\mR{{\mathbb R}}

\def\bP{{\mathbf P}}

\def\1{{\mathbf{1}}}

\def\sD{{\mathscr D}}
\def\sE{{\mathscr E}}
\def\sF{{\mathscr F}}
\def\sG{{\mathscr G}}

\def\sL{{\mathscr L}}

\def\sP{{\mathscr P}}

\def\geq{\geqslant}
\def\leq{\leqslant}

\def\eps{\varepsilon}

\def\e{\mathrm{e}}

\def\v{\mathrm{v}}
\def\u{\mathbf{u}}

\def\p{\partial}

\def\[{{\Big[}}
\def\]{{\Big]}}
\def\<{{\langle}}
\def\>{{\rangle}}
\def\({{\Big(}}
\def\){{\Big)}}

\def\bx{{\mathbf{x}}}
\def\tr{\mathrm {tr}}

\def\dif{{\mathord{{\rm d}}}}

\def\no{\nonumber}
\def\={&\!\!=\!\!&}
\def\bt{\begin{theorem}}
\def\et{\end{theorem}}
\def\bl{\begin{lemma}}
\def\el{\end{lemma}}
\def\br{\begin{remark}}
\def\er{\end{remark}}
\def\bx{\begin{Examples}}
\def\ex{\end{Examples}}
\def\bd{\begin{definition}}
\def\ed{\end{definition}}
\def\bp{\begin{proposition}}
\def\ep{\end{proposition}}
\def\bc{\begin{corollary}}
\def\ec{\end{corollary}}

\def\geq{\geqslant}
\def\leq{\leqslant}

\def\bP{{\mathbf P}}

\def\<{\langle} \def\>{\rangle}

 \def\beq{\begin{equation}}  
 
\def\e{\text{\rm{e}}}

\allowdisplaybreaks

\begin{document}

\title{Stochastic Hamiltonian flows with singular coefficients}

\author{XICHENG ZHANG
\\
\\
{\it Dedicated to the 60th birthday of Professor Michael R\"ockner} }

\thanks{This work is supported by NNSFs of China (Nos. 11271294, 11325105)}
\date{}

\address{
School of Mathematics and Statistics,
Wuhan University, Wuhan, Hubei 430072, P.R.China\\
Email: XichengZhang@gmail.com
 }

\begin{abstract}
In this paper we study the following stochastic Hamiltonian system in $\mR^{2d}$ (a second order stochastic differential equation),
$$
\dif \dot X_t=b(X_t,\dot X_t)\dif t+\sigma(X_t,\dot X_t)\dif W_t,\ \ (X_0,\dot X_0)=(x,\v)\in\mR^{2d},
$$
where $b(x,\v):\mR^{2d}\to\mR^d$ and $\sigma(x,\v):\mR^{2d}\to\mR^d\otimes\mR^d$ are two Borel measurable functions.
We show that if $\sigma$ is bounded and uniformly non-degenerate, and $b\in H^{2/3,0}_p$ and $\nabla\sigma\in L^p$ for some $p>2(2d+1)$, 
where $H^{\alpha,\beta}_p$ is the Bessel potential space with differentiability indices $\alpha$ in $x$ and $\beta$ in $\v$,
then the above stochastic equation  admits a unique strong solution so that $(x,\v)\mapsto Z_t(x,\v):=(X_t,\dot X_t)(x,\v)$ 
forms a stochastic homeomorphism flow, and $(x,\v)\mapsto Z_t(x,\v)$ is weakly differentiable with
ess.$\sup_{x,\v}\mE\left(\sup_{t\in[0,T]}|\nabla Z_t(x,\v)|^q\right)<\infty$ for all $q\geq 1$ and $T\geq 0$. Moreover, we also  show the uniqueness of 
probability measure-valued solutions for kinetic Fokker-Planck equations with rough coefficients 
by showing the well-posedness of the associated martingale problem and using the superposition principle established by Figalli \cite{Fi} and Trevisan \cite{Tre}.

\end{abstract}

\keywords{Stochastic Hamiltonian system, Weak differentiability, Krylov's estimate, Zvonkin's transformation, Kinetic Fokker-Planck operator.}
\maketitle \rm

\section{Introduction}
Consider the following {\it second order} time dependent stochastic differential equation (abbreviated as SDE):
\begin{align*}
\dif \dot X_t=b_t(X_t,\dot X_t)\dif t+\sigma_t(X_t,\dot X_t)\dif W_t,\ \ (X_0,\dot X_0)=(x,\v)\in\mR^{2d},
\end{align*}
where $b_t(x,\v):\mR_+\times\mR^{2d}\to\mR^d$ and $\sigma_t(x,\v):\mR_+\times\mR^{2d}\to\mR^d\otimes\mR^d$ are two Borel measurable functions,  $\dot X_t$ denotes the
first order derivative of $X_t$ with respect to $t$, and $W_t$ is a $d$-dimensional standard Brownian motion on some probability space $(\Omega,\sF,\mP)$. 
When $\sigma=0$, the above equation is the classical Newtonian mechanic equation, 
which describes the motion of a particle. When $\sigma\not=0$, 
it means that the motion is perturbed by some random external force. More backgrounds about the above stochastic Hamiltonian system are referred to \cite{So, Ta}, etc.

It is noticed that if we let $Z_t:=(X_t,\dot X_t)$, then $Z_t$ solves the following {\it one order} (degenerate) SDE:
\begin{align}\label{SDE0}
\dif Z_t=(\dot X_t,b_t(Z_t))\dif t+(0,\sigma_t(Z_t)\dif W_t),\ \ Z_0=z=(x,\v)\in\mR^{2d},
\end{align}
and whose time-dependent infinitesimal generator is given by
\begin{align}\label{LL}
\sL^{a,b}_t f(x,\v):=\tr(a_t\cdot\nabla^2_\v f)(x,\v)+(\v\cdot\nabla_x f)(x,\v)+(b_t\cdot\nabla_\v f)(x,\v).
\end{align}
Here $a_t(x,\v):=\tfrac{1}{2}(\sigma_t\sigma^*_t)(x,\v)$, $\nabla^2_\v f(x,\v)$ stands for the Hessian matrix, the asterisk and tr($\cdot$) denote 
the transpose and the trace of a matrix respectively.
Moreover, let $\mu_t$ be the probability distributional measure of $Z_t$ in $\mR^{2d}$. 
By It\^o's formula, one knows that $\mu_t$ solves the following Fokker-Planck equation in
the distributional sense:
\begin{align}\label{FP}
\p_t\mu_t=(\sL^{a,b}_t)^*\mu_t,\ \ \mu_0=\delta_z,
\end{align}
where $\delta_z$ is the Dirac measure at $z$. More precisely, for any $f\in C^2_c(\mR^{2d})$,
$$
\p_t\mu_t(f)=\mu_t(\sL^{a,b}_t f),\ \mu_0(f)=f(z),
$$
where $\mu_t(f)=\int f\dif\mu_t=\mE f(Z_t)$.
In the literature $\sL^{a,b}_t$ is also called kinetic Fokker-Planck or Kolmogorov's operator. 

During the past decade, there is an increasing interest in the study of  SDEs with singular or rough coefficients. 
In the non-degenerate case, Krylov and R\"ockner \cite{Kr-Ro}
showed the strong uniqueness to the following SDE in $\mR^d$:
$$
\dif X_t=b_t(X_t)\dif t+\dif W_t, \ X_0=x,
$$
where $b\in L^q_{lock}(\mR^+;L^p(\mR^d))$ with $\frac{d}{p}+\frac{2}{q}<1$. 
The argument in \cite{Kr-Ro} is based on Girsanov's theorem and some estimates from
the theory of PDE. In this framework, Fedrizzi and Flandoli  \cite{Fe-Fl, Fe-Fl0} studied the well-posedness of stochastic transport equations with rough coefficients. 
When $b$ is bounded measurable, the Malliavin differentiability of $X_t$ with respect to sample path $\omega$ and the weak differentiability of $X_t$
with respect to starting point $x$ were recently  studied in  \cite{Me} and \cite{Mo-Ni-Pr} respectively.
We also mention that weak uniqueness was studied in \cite{Ba-Ch} and \cite{Ji} under rather weak assumptions on $b$ (belonging to some Kato's class). 
Moreover, the multiplicative noise case was studied in \cite{Zh0, Zh2, Zh4} by using Zvonkin's transformation \cite{Zv} 
and some careful estimates of second order parabolic equations.

In the degenerate case, Chaudru de Raynal \cite{Ch} firstly showed the strong well-posedness for SDE \eqref{SDE0} under the assumptions that $\sigma$ is 
Lipschitz continuous and $b$ is $\alpha$-H\"older continuous in $x$ and $\beta$-H\"older continuous in $\v$ with $\alpha\in (\frac{2}{3},1)$ and $\beta\in(0,1)$.
The proofs in \cite{Ch} strongly depend on some explicit estimates for Kolmogorov operator with constant coefficients and Zvonkin's transformation.
In a recent joint work \cite{Wa-Zh} with F.Y. Wang, we also showed the strong uniqueness and homeomorphism property for \eqref{SDE0} 
under weaker H\"older-Dini's continuity assumption on $b$. The proofs in \cite{Wa-Zh} rely on a characterization of H\"older-Dini's spaces  and gradient 
estimates for the semigroup associated with the kinetic operator. Notice that in \cite{Ch} and \cite{Wa-Zh}, more general degenerate SDEs are considered,
while, the case with critical differentiability indices $\alpha=\frac{2}{3}$ and $\beta=0$ is left open. 

The purpose of this work is to establish a similar theory for degenerate SDE \eqref{SDE0} as in Krylov and R\"ockner's paper \cite{Kr-Ro} (see also \cite{Zh4}).
In particular, the critical indices $\alpha=\frac{2}{3}$ and $\beta=0$ are covered.
More precisely, we aim to prove that
\bt\label{Main}
Suppose that for some $K\geq 1$ and all $(t,x,\v)\in\mR_+\times\mR^{2d}$,
\begin{align*}
K^{-1}|\xi|\leq |\sigma^*_t(x,\v)\xi|\leq K|\xi|,\ \  \forall\xi\in\mR^d,\tag{{\bf UE}}
\end{align*}
where $\sigma^*$ denotes the transpose of matrix $\sigma$, and for some $p>2(2d+1)$,
$$
\kappa_0:=\sup_{s\geq 0}\|\nabla\sigma_s\|^p_p+\int_0^\infty\|(\mI-\Delta_x)^{\frac{1}{3}} b_s\|_p^p\dif s<\infty.
$$
Then for any $z=(x,\v)\in\mR^{2d}$, SDE \eqref{SDE0} admits a unique strong solution $Z_t(z)=(X_t,\dot X_t)$ so that $(t,z)\mapsto Z_t(z)$ has a bi-continuous version.
Moreover,
\begin{enumerate}[{\bf (A)}]
\item There is a null set $N$ such that for all $\omega\notin N$ and for each $t\geq 0$,
the map $z\mapsto Z_t(z,\omega)$ is a homeomorphism on $\mR^{2d}$.
\item For each $t\geq 0$, the map $z\mapsto Z_t(z)$ is weakly differentiable a.s., and for any $q\geq 1$ and $T>0$,
\begin{align}\label{YN1}
\mathrm{ess.}\sup_{z}\mE\left(\sup_{t\in[0,T]}|\nabla Z_t(z)|^q\right)<\infty,
\end{align}
where $\nabla$ denotes the generalized gradient.
\item Let $\sigma^n$ and $b^n$ be the regularized approximations of $\sigma$ and $b$ (see \eqref{Reg} below for definitions). 
Let $Z^n$ be the corresponding solution of SDE \eqref{SDE0} associated with $(\sigma^n,b^n)$. 
For any $q\geq 1$ and $T>0$, there exits a constant $C>0$ only depending on $T,K,\kappa_0,d,p,q$ 
such that
$$
\mE\left(\sup_{t\in[0,T]}|Z^n_t-Z_t|^q\right)\leq C\left(\|b^n-b\|^q_{\mL^p(T)}+n^{(\frac{2d}{p}-1)q}\right),\ \ n\in\mN.
$$
\end{enumerate}
\et

As a corollary, we have the following local well-posedness result by a standard localization argument.
\bc
Suppose that for any $T, R>0$, there exists a constant $K_{T,R}\geq 1$  such that for all $(t,x,\v)\in[0,T]\times B_R$ and $\xi\in\mR^d$,
\begin{align}
K_{T,R}^{-1}|\xi|\leq |\sigma^*_t(x,\v)\xi|\leq K_{T,R}|\xi|,\label{Con22}
\end{align}
where $B_R:=\{(x,\v): |(x,\v)|\leq R\}$,  and for some $p>2(2d+1)$,
$$
\sup_{t\in[0,T]}\|\nabla(\sigma_t\chi_R)\|^p_{p}+\int_0^T\|(\mI-\Delta_x)^{\frac{1}{3}} (b_s\chi_R)\|_p^p\dif s\leq K_{T,R},
$$
where $\chi_R:\mR^{2d}\to[0,1]$ is a smooth function with $\chi_R(z)=1$ for $|z|\leq R$ and $\chi_R(z)=0$ for $|z|>2R$.
Then for any fixed $(x,\v)\in\mR^{2d}$, SDE \eqref{SDE0} admits a unique local strong solution $(X_t,\dot X_t)$ up to the explosion time $\zeta$.
\ec
\begin{proof}
Let 
$$
\sigma^R_t(z):=\sigma_t(z\chi_R(z)),\ \ b^R_t(z):=b_t(z)\chi_R(z).
$$
By the assumptions, one sees that $(\sigma^R, b^R)$ satisfies the conditions of Theorem \ref{Main}. Hence, there exists a unique solution to 
the following SDE:
$$
\dif Z^R_t=(\dot X^R_t,b^R_t(Z^R_t))\dif t+(0,\sigma^R_t(Z^R_t)\dif W_t),\ \ Z^R_0=z=(x,\v)\in\mR^{2d},
$$
where $Z^R_t=(X^R_t,\dot X^R_t)$.
Define
$$
\zeta_R:=\inf\Big\{t\geq 0: |Z^R_t|\geq R\Big\},\ \ Z_t:=Z^R_t,\ t\in[0,\zeta_R].
$$
Since $Z^{R'}_t|_{[0,\zeta_R]}=Z^R_t|_{[0,\zeta_R]}$ for $R'>R$, one sees that $R\mapsto \zeta_R$ is increasing and the above $Z_t$ is well-defined. 
Clearly, $\zeta=\lim_{R\to\infty}\zeta_R$ is the explosion time of $Z_t$, and $Z_t$ uniquely solves \eqref{SDE0} before $\zeta$.
\end{proof}

The strategy of proving Theorem \ref{Main} is still based on Zvonkin's transformation. As in the non-degenerate case \cite{Zh4}, 
we need to establish the $L^p$-maximal regularity estimate  to the following degenerate 
parabolic equation (see Theorem \ref{Th32} below):
$$
\p_t u=\sL^{a,b}_t u+f,\ \ u_0=0.
$$ 
Here we shall use the freezing coefficient argument and the $L^p$-estimate established in \cite{Br-Cu-La-Pr} and \cite{Bo} 
for degenerate operators with  constant coefficients (see also \cite{Ch-Zh} for the case of nonlocal operators). 
Compared with \cite{Br-Cu-La-Pr} and \cite{Pr}, we not only consider the optimal regularity of $u$ along the nondegenerate $\v$-direction,
but also the optimal regularity of $u$ along the degenerate $x$-direction.


On the other hand, from the viewpoint of PDE, the well-posedness of Fokker-Planck equation \eqref{FP} (especially uniqueness) with rough coefficients 
is a quite involved problem. Since $a$ and $b$ possess less regularities and $\sL^{a,b}_t$ is a degenerate operator, 
the direct analytical approach seems not work (cf. \cite{Bo-Kr-Ro1, Bo-Kr-Ro2}). 
Let $\cP(\mR^{2d})$ be the set of all probability measures on $\mR^{2d}$. We shall use a probabilistic method
to prove the following result.
\bt\label{Main2}
Suppose that $\sigma$ satisfies {\bf (UE)} and for any $T>0$,
$$
\lim_{|z-z'|\to 0}\sup_{t\in[0,T]}\|\sigma_t(z)-\sigma_t(z')\|=0,
$$
and $b\in L^q_{loc}(\mR_+; L^q(\mR^{2d}))$ for some $q\in(2(2d+1),\infty]$. 
Then for any $\nu\in\cP(\mR^{2d})$, there exists a unique probability measure-valued solution 
$\mu_t\in\cP(\mR^{2d})$ to \eqref{FP} in the distributional sense in the class that
 $t\mapsto\mu_t$ is weakly continuous with $\mu_0=\nu$ and
 $$
 \int^t_0\!\!\int_{\mR^{2d}}(|\v|+|b_s(x,\v)|)\mu_s(\dif x,\dif\v)\dif s<\infty,\ \ t>0.
 $$
\et

The proof of this result is based on Figalli and Trevisan's superposition characterization for the solutions of 
Fokker-Planck equation in terms of martingale problem associated with $\sigma$ and $b$. 
More precisely, Figalli \cite{Fi} and Trevisan \cite{Tre} showed that for any weakly continuous probability measure-valued 
solution $\mu_t$ of \eqref{FP} with initial value $\nu\in\cP(\mR^{2d})$, 
there exists a martingale solution for operator $\sL^{a,b}_t$ (a probability measure $\mP_\nu$ over the space of 
all continuous functions from $\mR_+$ to $\mR^{2d}$ denoted by $\Omega$) such that for all $t\in\mR_+$,
$$
\int_{\mR^{2d}}\varphi(z)\mu_t(\dif z)=\int_{\Omega}\varphi(\omega_t)\mP_\nu(\dif \omega),
$$
where $t\mapsto\omega_t$ is the coordinate process over $\Omega$.
Hence, in order to prove Theorem \ref{Main2}, it suffices to show the well-posedness of martingale problem for $\sL^{a,b}_t$ 
in the sense of Stroock and Varadhan \cite{St-Va}.
This will be achieved by proving some Krylov's type estimate (see Theorem \ref{Th41} below), which is also a key tool for proving Theorem \ref{Main}.
It is remarked that in \cite{Ro-Zh}, we have already used this technique to show the uniqueness of measure-valued solutions and $L^p$-solutions to 
possibly degenerate second order Fokker-Planck equations under some weak conditions on the coefficients (but not the case of Theorem \ref{Main2}).

This paper is organized as follows: in Section 2, we introduce some anisotropic fractional Bessel potential spaces,
and prepare some useful estimates for later use. In Section 3, we show the $L^p$-maximal regularity estimate
for kinetic Fokker-Planck equations. In Section 4, we study the martingale problem associated with $(\sigma,b)$ 
under the same assumptions as in Theorem \ref{Main2} by showing the basic Krylov's type estimate.
In particular, we first prove Theorem \ref{Main2}.
In Section 5, we then prove Theorem \ref{Main} by
using Zvonkin's transformation and Krylov's estimate obtained in the previous section. 
In Appendix, a stochastic Gronwall's type lemma used in Section 5 is given.

Convention: The letter $C$ with or without subscripts will denote an unimportant constant, whose value may change in different places. Moreover, $A\preceq B$
means that $A\leq CB$ for some constant $C>0$, and $A\asymp B$ means that  $C^{-1}B\leq A\leq CB$ for some $C>1$.

After this work was finished, I was informed by Professor Enrico Priola during ``The 8th International Conference on Stochastic Analysis and its Applications" 
held at BIT that, very recently, Fedrizzi, Flandoli, Priola and Vovelle \cite{Fe-Fl-Pr-Vo} also obtained 
the strong well-posedness together with their flow property of SDE \eqref{SDE0} under the conditions $\sigma_t(z)=\mI$ and $b_t(z)=b(z)$ possessing the following regularity
$$
\|(\mI-\Delta_x)^{s/2}b\|_p<\infty
$$
for some $s>2/3$ and $p>6d$. 
\section{Preliminaries}

For $\alpha\geq 0$ and $p\in(1,\infty)$, let $H^\alpha_p:=H^\alpha_p(\mR^d):=(\mI-\Delta)^{-\frac{\alpha}{2}}(L^p(\mR^d))$ 
be the usual Bessel potential space with norm
$$
\|f\|_{\alpha,p}:=\|(\mI-\Delta)^{\frac{\alpha}{2}} f\|_p,
$$
where $\|\cdot\|_p$ is the usual $L^p$-norm, and $\Delta$ is the Laplacian. For $\alpha\in(0,2)$, let $\Delta^{\frac{\alpha}{2}}:=-(-\Delta)^{\frac{\alpha}{2}}$ 
be the usual fractional Laplacian. Notice that up to a constant $C(\alpha,d)>0$, an alternative definition of  $\Delta^{\frac{\alpha}{2}}$ is given by
\begin{align}\label{Def1}
\Delta^{\frac{\alpha}{2}} f(x):=\lim_{\eps\downarrow 0}\int_{|y|\geq\eps}\delta_yf(x)|y|^{-d-\alpha}\dif y,\ \ \delta_yf(x):=f(x+y)-f(x).
\end{align}
We will frequently use such a definition below.
It is well-known that by the boundedness of Riesz's transformation (cf. \cite{St}), 
\begin{align}\label{Riesz}
\|\Delta^{\frac{1}{2}}f\|_p\asymp\|\nabla f\|_p,\ p>1,
\end{align}
and an equivalent norm in $H^\alpha_p$ is given by
\begin{align}\label{Equiv}
\|f\|_{\alpha,p}\asymp\|f\|_{p}+\|\Delta^{\frac{\alpha-[\alpha]}{2}}\nabla^{[\alpha]} f\|_p,
\end{align}
where $[\alpha]$ is the integer part of real number $\alpha$, and we have used the convention $\Delta^0:=\mI$.
Notice that for $\alpha\in(0,1]$ and $p>1$, 
\begin{align}\label{LI0}
\|f(\cdot+x)-f(\cdot)\|_p\preceq \|\Delta^{\frac{\alpha}{2}}f\|_p|x|^\alpha,
\end{align}
and in particular,
\begin{align}\label{LI}
\|f(\cdot+x)-f(\cdot)\|_p\preceq \|f\|_{\alpha,p}(|x|^\alpha\wedge 1).
\end{align}
Moreover, we also have the following interpolation inequality: for any $0\leq\alpha<\beta<\infty$,
\begin{align}\label{Inter}
\|f\|_{\alpha,p}\leq C(p,d,\alpha,\beta)\|f\|^{\frac{\beta-\alpha}{\beta}}_p\|f\|^{\frac{\alpha}{\beta}}_{\beta,p},
\end{align}
and the following Sobolev embedding results hold: for any $\alpha\in(0,1)$, if $p\alpha>d$, then
\begin{align}\label{Em0}
\|f\|_\infty+\sup_{x\not=y}\frac{|f(x)-f(y)|}{|x-y|^\gamma}\leq C(p,d,\alpha,\gamma)\|f\|_{\alpha,p},\ \gamma\in(0,\alpha-\tfrac{d}{p});
\end{align}
if $p\alpha<d$, then
\begin{align}\label{Em1}
\|f\|_{q}\leq C(p,d,\alpha,q)\|f\|_{\alpha,p},\ \ q\in[p,\tfrac{pd}{d-p\alpha}].
\end{align}
All the above facts are standard and can be found in \cite{Be-Lo} and \cite{St}.

\

To treat the kinetic Fokker-Planck operator, we introduce the following anisotropic Bessel potential spaces.
Let $C^\infty_c(\mR^{2d})$ be the space of all smooth functions on $\mR^{2d}$ with compact supports. 
For $\alpha,\beta\geq 0$, we define the Bessel potential space
$H^{\alpha,\beta}_p:=H^{\alpha,\beta}_p(\mR^{2d})$ as the completion of $C^\infty_c(\mR^{2d})$ with respect to norm:
$$
\|f\|_{\alpha,\beta; p}:=\|(\mI-\Delta_x)^{\frac{\alpha}{2}}f\|_p+\|(\mI-\Delta_\v)^{\frac{\beta}{2}}f\|_p.
$$
Notice that by the Mihlin multiplier theorem (cf. \cite{Be-Lo}),
$$
\|f\|_{\alpha,\beta;p}\asymp \|f\|_p+\|\Delta_x^{\frac{\alpha}{2}}f\|_p+\|\Delta_\v^{\frac{\beta}{2}}f\|_p
\asymp \|((\mI-\Delta_x)^{\frac{\alpha}{2}}+(\mI-\Delta_\v)^{\frac{\beta}{2}})f\|_p.
$$
In the following, we simply write
$$
H^{\infty,\infty}_p:=H^{\infty,\infty}_p(\mR^{2d}):=\cap_{\alpha,\beta\geq 0}H^{\alpha,\beta}_p(\mR^{2d}).
$$
\bl
(i) For any $\alpha,\beta\geq 0$, $\theta\in[0,1]$ and $p>1$, there is a constant $C=C(\alpha,\beta,\theta,p,d)>0$ such that
\begin{align}
&\|(\mI-\Delta_x)^{\frac{\theta\alpha}{2}}(\mI-\Delta_\v)^{\frac{(1-\theta)\beta}{2}}f\|_p\leq C\|f\|_{\alpha,\beta;p},\label{Em8}\\
&\quad \|\Delta_x^{\frac{\theta\alpha}{2}}\Delta_\v^{\frac{(1-\theta)\alpha}{2}}f\|_p\leq C\|(\Delta_x+\Delta_\v)^{\frac{\alpha}{2}}f\|_p.\label{Em88}
\end{align}
In particular, for any $\alpha\geq 0$, $\beta\geq 1$ and $p>1$, we have
\begin{align}\label{E28}
\|\nabla_\v f\|_{\alpha(\beta-1)/\beta,\beta-1;p}\leq C\|f\|_{\alpha,\beta;p}.
\end{align}
(ii) Let $\alpha,\beta\geq 0$ and $p>1$ with $d\not=\frac{p\alpha\beta}{\alpha+\beta}$. Set 
$$
p^*:=
\left\{
\begin{aligned}
&dp/(d-\tfrac{p\alpha\beta}{\alpha+\beta}),&\ \ d>\tfrac{p\alpha\beta}{\alpha+\beta};\\
&\infty,&\ \ d<\tfrac{p\alpha\beta}{\alpha+\beta}. 
\end{aligned}
\right.
$$
For any $q\in[p,p^*]$,  there is a constant $C=C(\alpha,\beta,p,q,d)>0$ such that
\begin{align}\label{E29}
\|f\|_{q}\leq C\|f\|_{\alpha,\beta;p}.
\end{align} 
\el
\begin{proof}
(i) It follows by the Mihlin multiplier theorem and \eqref{Riesz}.

(ii) For \eqref{E29}, by \eqref{Em0}, \eqref{Em1} and \eqref{Em8} with $\theta=\frac{\beta}{\alpha+\beta}$, we have
\begin{align*}
\|f\|^q_q&=\int_{\mR^d}\|f(\cdot,\v)\|_q^q\dif \v \preceq
\int_{\mR^d}\|(\mI-\Delta_x)^{\frac{\theta\alpha}{2}}f(\cdot,\v)\|_p^q\dif \v\\
&\leq\left(\int_{\mR^d}\left(\int_{\mR^d}|(\mI-\Delta_x)^{\frac{\theta\alpha}{2}}f(x,\v)|^q\dif\v\right)^{p/q}\dif x\right)^{q/p}\\
&\preceq\left(\int_{\mR^d}\|(\mI-\Delta_\v)^{\frac{(1-\theta)\beta}{2}}(\mI-\Delta_x)^{\frac{\theta\alpha}{2}}f(x,\cdot)\|^p_{p}\dif x\right)^{q/p}
\preceq\|f\|^q_{\alpha,\beta;p},
\end{align*}
where the second inequality is due to Minkovskii's inequality.
\end{proof}
Let $a:\mR^{2d}\to\mR^d\otimes\mR^d$ be a measurable function. Write
$$
\sL^a u:=\mathrm{tr}(a\cdot\nabla^2_\v u)+\v\cdot\nabla_x u.
$$
We have
\bl\label{Le22}
Let $\alpha\in(0,1)$ and $p>d/\alpha$. Suppose that
$$
\kappa_0:=\sup_{\v}\|\Delta_x^{\frac{\alpha}{2}}a(\cdot,\v)\|_p+\|a\|_\infty<\infty.
$$ 
For any $\eps\in(0,1)$, there is a constant $C_\eps=C_\eps(p,d,\alpha,\kappa_0)>0$ such that for all $u\in H^{\infty,\infty}_p$,
\begin{align}\label{E210}
\|[\Delta_x^{\frac{\alpha}{2}},\sL^a]u\|_p\leq\eps\|\Delta_x^{\frac{\alpha}{2}}\nabla^2_\v u\|_{p}+C_\eps \|\nabla^2_\v u\|_{p},
\end{align}
where $[\Delta_x^{\frac{\alpha}{2}},\sL^a]u:=\Delta_x^{\frac{\alpha}{2}}(\sL^a u)-\sL^a(\Delta_x^{\frac{\alpha}{2}}u)$.
\el
\begin{proof}
Notice that by definition \eqref{Def1},
$$
[\Delta_x^{\frac{\alpha}{2}},\sL^a]u=\tr\big(\Delta_x^{\frac{\alpha}{2}} a\cdot\nabla^2_\v u\big)
+\int_{\mR^d}\tr\big(\delta_{(y,0)}a\cdot\nabla^2_\v\delta_{(y,0)}u\big)|y|^{-d-\alpha}\dif y.
$$
Hence,
$$
\|[\Delta_x^{\frac{\alpha}{2}},\sL^a]u\|_p\leq\|\tr(\Delta_x^{\frac{\alpha}{2}} a\cdot\nabla^2_\v u)\|_p
+\int_{\mR^d}\big\|\tr\big(\delta_{(y,0)}a\cdot\nabla^2_\v\delta_{(y,0)}u\big)\big\|_p|y|^{-d-\alpha}\dif y.
$$
Let $\beta\in(\frac{d}{p},\alpha)$. By \eqref{Em0}, we have
\begin{align*}
\|\tr(\Delta_x^{\frac{\alpha}{2}} a\cdot\nabla^2_\v u)\|_p^p
&\preceq\int_{\mR^d}\|\Delta_x^{\frac{\alpha}{2}} a(\cdot,\v)\|^p_{p}\|\nabla^2_\v u(\cdot,\v)\|_\infty^p\dif \v\\
&\preceq\sup_{\v}\|\Delta_x^{\frac{\alpha}{2}}a(\cdot,\v)\|^p_p\|\nabla^2_\v u\|_{\beta,0; p}^p,
\end{align*}
and for  $\gamma\in(0,\beta-\frac{d}{p})$,
\begin{align*}
&\big\|\tr\big(\delta_{(y,0)}a\cdot\nabla^2_\v\delta_{(y,0)}u\big)\big\|_p^p
\preceq\int_{\mR^d}\|\delta_y a(\cdot,\v)\|_{p}^p\cdot\|\nabla^2_\v\delta_y u(\cdot,\v)\|_{\infty}^p\dif \v\\
&\qquad\stackrel{\eqref{LI0}}{\preceq}|y|^{(\alpha+\gamma)p}\int_{\mR^d}\|\Delta_x^{\frac{\alpha}{2}}a(\cdot,\v)\|_p^p\|\nabla^2_\v u(\cdot,\v)\|_{\beta,p}^p\dif\v\\
&\qquad\ \preceq |y|^{(\alpha+\gamma)p}\sup_{\v}\|\Delta_x^{\frac{\alpha}{2}}a(\cdot,\v)\|^p_p\|\nabla^2_\v u\|_{\beta,0;p}^p.
\end{align*}
Moreover, it is easy to see that
$$
\big\|\tr\big(\delta_{(y,0)}a\cdot\nabla^2_\v\delta_{(y,0)}u\big)\big\|_p\preceq\|a\|_\infty\|\nabla^2_\v u\|_p.
$$
Therefore,
$$
\big\|\tr\big(\delta_{(y,0)}a\cdot\nabla^2_\v\delta_{(y,0)}u\big)\big\|_p\preceq \kappa_0(|y|^{\alpha+\gamma}\wedge 1)\|\nabla^2_\v u\|_{\beta,0;p}.
$$
Combining the above calculations, we get for some $C=C(p,d,\alpha,\beta)>0$,
\begin{align}\label{LJ1}
\|[\Delta_x^{\frac{\alpha}{2}},\sL^a]u\|_p\leq C\kappa_0\|\nabla^2_\v u\|_{\beta,0;p},
\end{align}
On the other hand, by the interpolation inequality \eqref{Inter} and Young's inequality, we have for any $\eps\in(0,1)$,
$$
\|\nabla^2_\v u\|_{\beta,0;p}\preceq \|\nabla^2_\v u\|^\frac{\beta}{\alpha}_{\alpha,0;p}\|\nabla^2_\v u\|^\frac{\alpha-\beta}{\alpha}_{p}
\leq\eps\|\nabla^2_\v u\|_{\alpha,0;p}+C_\eps\|\nabla^2_\v u\|_{p}.
$$
Estimate \eqref{E210} now follows by \eqref{LJ1}.
\end{proof}
\bl\label{Le23}
For any $\alpha,\beta\in(0,1)$ and $p>(\alpha+\beta)d/(\alpha\beta)$,
there is a constant $C=C(\alpha,\beta,p,d)>0$ such that for all $b\in H^{\alpha,0}_p$ and $u\in H^{\infty,\infty}_p$,
$$
\|b\cdot\nabla_\v u\|_{\alpha,0; p}\leq C\|b\|_{\alpha,0;p}\Big(\|\Delta_x^{\frac{\alpha}{2}}\nabla_\v u\|_{0,\beta;p}+\|\nabla_\v u\|_{0,\beta; p}\Big).
$$
\el
\begin{proof}
Notice that by \eqref{Equiv},
\begin{align*}
\|b\cdot\nabla_\v u\|_{\alpha,0;p}\preceq\|b\cdot\nabla_\v u\|_p+\|\Delta_x^{\frac{\alpha}{2}}(b\cdot\nabla_\v u)\|_p.
\end{align*}
By definition \eqref{Def1}, we have
\begin{align*}
&\|\Delta_x^{\frac{\alpha}{2}}(b\cdot\nabla_\v u)\|_p\leq\|(\Delta_x^{\frac{\alpha}{2}} b)\cdot\nabla_\v u\|_p
+\|b\cdot\nabla_\v \Delta_x^{\frac{\alpha}{2}}u\|_p\\
&\qquad+\int_{\mR^d}\big\|\delta_{(y,0)}b\cdot\nabla_\v\delta_{(y,0)}u\big\|_p|y|^{-d-\alpha}\dif y=:I_1+I_2+I_3.
\end{align*}
For $I_1$, since $p>(\alpha+\beta)d/(\alpha\beta)$, by \eqref{E29} with $q=\infty$, we have
\begin{align*}
I_1\preceq \|\Delta_x^{\frac{\alpha}{2}} b\|_p\|\nabla_\v u\|_\infty\preceq\|b\|_{\alpha,0; p}\|\nabla_\v u\|_{\alpha,\beta;p}.
\end{align*}
For $I_2$, since $p\alpha>d$, by \eqref{Em0} we have
\begin{align*}
I_2^p&=\int_{\mR^d}\|b(\cdot,\v)\cdot\nabla_\v\Delta_x^{\frac{\alpha}{2}} u(\cdot,\v)\|_p^p\dif \v\\
&\preceq\int_{\mR^d}\|b(\cdot,\v)\|_\infty^p\|\nabla_\v\Delta_x^{\frac{\alpha}{2}} u(\cdot,\v)\|_p^p\dif \v\\
&\preceq\int_{\mR^d}\|b(\cdot,\v)\|_{\alpha,p}^p\dif \v\sup_\v\|\nabla_\v\Delta_x^{\frac{\alpha}{2}} u(\cdot,\v)\|_p^p.
\end{align*}
For $I_3$, by \eqref{LI} and \eqref{Em0} again, we have for any $\gamma\in(0,\alpha-\frac{d}{p})$, 
\begin{align*}
&\int_{\mR^d}\|\delta_yb(\cdot,\v)\cdot\nabla_\v\delta_yu(\cdot,\v)\|_p^p\dif \v
\preceq\int_{\mR^d}\|\delta_y b(\cdot,\v)\|_p^p\|\delta_y\nabla_\v u(\cdot,\v)\|_\infty^p\dif \v\\
&\qquad\preceq\int_{\mR^d}\|b(\cdot,\v)\|_{\alpha,p}^p(|y|^{\alpha p}\wedge 1)\|\nabla_\v u(\cdot,\v)\|_{\alpha,p}^p(|y|^{\gamma p}\wedge 1)\dif \v\\
&\qquad\preceq\left(\int_{\mR^d}\|b(\cdot,\v)\|_{\alpha,p}^p\dif \v\right)\sup_\v\|\nabla_\v u(\cdot,\v)\|_{\alpha,p}^p(|y|^{(\alpha+\gamma)p}\wedge 1).
\end{align*}
On the other hand, notice that by $p\beta>d$ and \eqref{Em0}, 
\begin{align*}
\sup_\v\|\nabla_\v\Delta_x^{\frac{\alpha}{2}} u(\cdot,\v)\|_p^p
\leq\int_{\mR^d}\sup_\v |\nabla_\v\Delta_x^{\frac{\alpha}{2}} u(x,\v)|^p\dif x
\preceq\int_{\mR^d}\|\nabla_\v\Delta_x^{\frac{\alpha}{2}} u(x,\cdot)\|^p_{\beta, p}\dif x,
\end{align*}
and similarly,
$$
\|b\cdot\nabla_\v u\|^p_p\preceq\int_{\mR^d} \|b(\cdot,\v)\|^p_{\alpha,p}\|\nabla_\v u(\cdot,\v)\|^p_{p}\dif\v\preceq\|b\|^p_{\alpha,0;p}\|\nabla_\v u\|^p_{0,\beta;p}.
$$
Combining the above calculations, we obtain the desired estimate.
\end{proof}

Let $\varrho:\mR^{2d}\to[0,\infty)$ be a smooth function with support in the unit ball and $\int\varrho=1$. Define
\begin{align}\label{Rho}
\varrho_\eps(z):=\eps^{-2d}\varrho(\eps^{-1}z),\ \ \eps\in(0,1),
\end{align}
and for a locally integrable function $u:\mR^{2d}\to\mR$,
$$
u_{\eps}(z):=u*\varrho_\eps(z)=\int_{\mR^{2d}}u(z')\varrho_\eps(z-z')\dif z'.
$$
Let $\sP$ be an operator on the space of locally integrable functions. We define
\begin{align}\label{Def0}
[\varrho_\eps, \sP]u:=(\sP u)*\varrho_\eps-\sP(u*\varrho_\eps).
\end{align}
We need the following commutator estimate results.
\bl\label{Le24}
\begin{enumerate}[(i)]
\item Let $p\in[1,\infty)$ and $q,r\in[p,\infty]$ with $\frac{1}{p}=\frac{1}{q}+\frac{1}{r}$.
For any $b\in L^q(\mR^{2d})$ and $u\in H^{0,1}_r$, we have
\begin{align}\label{ES9}
\lim_{\eps\to 0}\|[\varrho_\eps, b\cdot\nabla_\v]u\|_p=0.
\end{align}
\item 
Let $a:\mR^{2d}\to\mR^d\otimes\mR^d$ be a bounded measurable function.
For any $p\in[1,\infty)$ and $u\in H^{0,2}_p$, we have
\begin{align}\label{ES8}
\lim_{\eps\to 0}\|[\varrho_\eps, \sL^a]u\|_p=0.
\end{align}
\end{enumerate}
\el
\begin{proof}
(i) It follows by \cite[Lemma 4.2]{Zh22}.

\

(ii) By definition, we can write for $z=(x,\v)$,
\begin{align*}
[\varrho_\eps, \sL^a]u(z)&=(\sL^au)*\varrho_\eps(z)-\sL^a(u*\varrho_\eps)(z)\\
&=\int_{\mR^{2d}}\tr((a(z')-a(z))\cdot\nabla^2_\v u(z'))\varrho_\eps(z-z')\dif z'\\
&\quad+\int_{\mR^{2d}}(\v'-\v) \cdot\nabla_x u(z')\varrho_\eps(z-z')\dif z'\\
&=:I^\eps_1(t,z)+I^\eps_2(t,z).
\end{align*}
For $I^\eps_1(t,z)$, by Jensen's inequality and the assumption, we have
\begin{align*}
|I^\eps_1(t,z)|^p&\leq \int_{\mR^{2d}}|\tr((a(z')-a(z))\cdot\nabla^2_\v u(z'))|^p\varrho_\eps(z-z')\dif z'\\
&\leq(2\|a\|_\infty)^p\int_{\mR^{2d}}\|\nabla^2_\v u(z')\|^p\varrho_\eps(z-z')\dif z'.
\end{align*}
For $I^\eps_2(t,z)$, by the integration by parts and H\"older's inequality, we have
\begin{align*}
|I^\eps_2(t,z)|^p&=\left|\int_{\mR^{2d}}(\v'-\v) \cdot \nabla_x\varrho_\eps(z-z')u(z')\dif z'\right|^p\\
&\leq \eps^{p}\left(\int_{\mR^{2d}}|\nabla_x\varrho_\eps(z-z')|\cdot|u(z')|\dif z'\right)^p\\
&\leq \eps^{p}\left(\int_{\mR^{2d}}|\nabla_x\varrho_\eps(z)|\dif z\right)^{p-1}
\int_{\mR^{2d}}|\nabla_x\varrho_\eps(z-z')|\cdot|u(z')|^p\dif z'\\
&=\left(\int_{\mR^{2d}}|\nabla_x\varrho(z)|\dif z\right)^{p-1}
\int_{\mR^{2d}}|(\nabla_x\varrho)_\eps(z-z')|\cdot|u(z')|^p\dif z'.
\end{align*}
Combining the above calculations, we get 
$$
\|[\varrho_\eps, \sL^a]u\|_p\leq C\|u\|_{0,2;p}.
$$
Hence, for any $u\in H^{0,2}_p$, it is easy to see that
\begin{align}\label{GH0}
\lim_{\eps'\to 0}\sup_{\eps\in(0,1)}\|[\varrho_\eps, \sL^a](u_{\eps'}-u)\|_p\leq C\lim_{\eps'\to 0}\|u_{\eps'}-u\|_{0,2;p}=0.
\end{align}
Moreover, for fixed $\eps'\in(0,1)$, since $u_{\eps'}\in H^{\infty,\infty}_p$, by \cite[Lemma 4.2]{Zh22} we have
$$
\lim_{\eps\to 0}\|[\varrho_\eps, \sL^a]u_{\eps'}\|_p=0,
$$
which together with \eqref{GH0} implies \eqref{ES8}.
\end{proof}

Let $\sigma_t(x,\v)=\sigma_t$ be independent of $(x,\v)$. Define for $t<s$,
\begin{align}\label{PP}
P_{t,s}f(x,\v)=\mE f(x+(s-t)\v+X_{t,s}, \v+V_{t,s}),
\end{align}
where
$$
(X_{t,s},V_{t,s})=\left(\int^s_tV_{t,r}\dif r, \int^s_t\sigma_r\dif W_r\right).
$$

We need the following basic $L^p$-regularity estimates related to $P_{t,s}$, which plays a basic role in the next section.
\bt\label{Th25}
Let $T>0$. Suppose that for some $K>0$ and all $t\in[0,T]$,  
$$
K^{-1}|\xi|\leq |\sigma^*_t\xi|\leq K|\xi|,\ \ \xi\in\mR^d.
$$
\begin{enumerate}[(i)]
\item For any $\alpha,\beta\geq 0$ and $p>1$, there exists a positive constant $C=C(K,T,p,d,\alpha,\beta)$ 
such that for all $f\in L^p(\mR^{2d})$ and $0\leq t<s\leq T$,
\begin{align}\label{Gr}
\begin{split}
\|P_{t,s}f\|_{\alpha,0;p}\leq C(s-t)^{-\frac{3\alpha}{2}}\|f\|_p,\\
\|P_{t,s}f\|_{0,\beta;p}\leq C(s-t)^{-\frac{\beta}{2}}\|f\|_p.
\end{split}
\end{align}
\item For any $p>1$, there exists a positive constant $C_p=C_p(K,d)$ such that for all $\lambda\geq 0$ and $f\in\mL^p(T)=L^p([0,T]\times\mR^{2d})$,
\begin{align}\label{Pri22}
\|\nabla^2_\v u^\lambda\|_{\mL^P(T)}+\|\Delta^{\frac{1}{3}}_x u^\lambda\|_{\mL^P(T)}\leq C_p\|f\|_{\mL^P(T)},
\end{align}
where $u^\lambda_t(x,\v):=\int^T_t\e^{\lambda(t-s)}P_{t,s}f_s(x,\v)\dif s$ satisfies 
$$
\p_t u^\lambda+\sL^{a,0}_tu^\lambda-\lambda u^\lambda+f=0
$$
in the distributional sense.
\end{enumerate}
\et
\begin{proof}
(i) It follows by the following gradient estimate and the interpolation theorem (see \cite[Theorem 2.10]{Wa-Zh}),
$$
\|\nabla^k_x\nabla^m_\v P_{t,s}f\|_p\leq C(s-t)^{-\frac{3k+m}{2}}\|f\|_p,\ k,m\in\mN_0.
$$
(ii) It is a consequence of \cite{Br-Cu-La-Pr} and \cite[Theorem 2.1]{Bo} (see also \cite[Theorem 3.3]{Ch-Zh}).
\end{proof}
\br
Notice that in the references \cite{Wa-Zh} and \cite{Ch-Zh}, the positions of $t$ and $s$ are exchanged.
\er

\section{Maximal $L^p$-solutions of kinetic Fokker-Planck equations}

Throughout this section, we fix $T>0$. 
Let $p\in(1,\infty)$ and $\alpha,\beta\geq 0$. For $t\in[0,T]$, we introduce the following Banach spaces with natural norms:
$$
\mL^p(t,T):=L^p([t,T];L^p(\mR^{2d})),\ \ \mH^{\alpha,\beta}_p(t,T):=L^p([t,T];H^{\alpha,\beta}_p(\mR^{2d})).
$$
For simplicity of notation,  we write
$$
\mL^p(T):=\mL^p(0,T),\ \ \mH^{\alpha,\beta}_p(T):=\mH^{\alpha,\beta}_p(0,T).
$$

We assume that $a:[0,T]\times\mR^{2d}\to\mR^d\otimes\mR^d$ is symmetric and satisfies
that for some $K\geq 1$ and $\delta\in(0,1)$,
\begin{align}\tag{\bf H$^{\delta,p}_K$}
\left\{
\begin{aligned}
&K^{-1}\cdot\mI\leq a_t(z)\leq K\cdot\mI,\,\,(t,z)\in[0,T]\times\mR^{2d}\\
&\omega_a(\delta):=\sup_{|z-z'|\leq\delta}\sup_{t\in[0,T]}\|a_t(z)-a_t(z')\|\leq\tfrac{1}{2(C_p+1)}
\end{aligned}
\right\},
\end{align}
where $C_p$ is the same as in \eqref{Pri22}.
Here and in the remainder of this paper, $\|\cdot\|$ denotes the Hilbert-Schmidt norm.
For $\lambda\geq 0$, consider the following backward kinetic Fokker-Planck equation
\begin{align}\label{PDE}
\p_t u+\sL^{a,b}_t u-\lambda u+f=0, \ \ u_T=0,
\end{align}
where $f_t(x,\v):[0,T]\times\mR^{2d}\to\mR$ is a Borel function. 
We first introduce the following notion of solutions to the above equation.
\bd\label{Def31}
Let $p\in(1,\infty)$ and $f\in\mL^p(T)$. A Borel function $u\in \mH^{0,2}_p(T)$ 
is called a solution of  \eqref{PDE} if for any $\varphi\in C^\infty_c(\mR^{2d})$ and all $t\in[0,T]$,
\begin{align}\label{Sol}
\begin{split}
\<u_t,\varphi\>&=\int^T_t\<\tr(a_s\cdot\nabla^2_\v u_s),\varphi\>\dif s-\int^T_t\<\v\cdot\nabla_x\varphi, u_s\>\dif s\\
&\quad+\int^T_t\<b_s\cdot\nabla_\v u_s,\varphi\>\dif s-\lambda\int^T_t\<u_s,\varphi\>\dif s+\int^T_t\<f_s,\varphi\>\dif s,
\end{split}
\end{align}
where $\<u_t,\varphi\>:=\int_{\mR^{2d}}u_t(z)\varphi(z)\dif z$.
\ed
The main aim of this section is to show that
\bt\label{Th32}
Let $\alpha\in[0,\frac{2}{3})$, $\beta\in(1,2)$ and $p>\frac{2}{(2-3\alpha)\wedge(2-\beta)}$ be not equal to $\frac{d(\alpha+\beta)}{\alpha(\beta-1)}$. 
Suppose that $a$ satisfies {\bf (H$^{\delta,p}_K$)}, and for some $q\in[p\vee \frac{d(\alpha+\beta)}{\alpha(\beta-1)},\infty]$,
$$
\kappa_0:=\|b\|_{L^p([0,T]; L^q(\mR^{2d}))}<\infty.
$$
\begin{enumerate}[(i)]
\item For any $f\in \mL^p(T)$, there exists a unique solution $u=u^\lambda$ to \eqref{PDE} in the sense of Definition \ref{Def31} with
\begin{align}\label{Pri0}
\|u^\lambda\|_{\mH^{2/3,2}_p(T)}\leq C\|f\|_{\mL^p(T)}, \ 
\end{align}
and for all $t\in[0,T]$,
\begin{align}\label{Pri00}
\|u^\lambda_t\|_{\alpha,\beta;p}\leq C((T-t)\wedge\lambda^{-1})^{1-\frac{1}{p}-\frac{(3\alpha)\vee\beta}{2}}\|f\|_{\mL^p(t,T)},
\end{align}
where the constant $C$ only depends on $d,\delta,K,\alpha,\beta,p, q,T$ and $\kappa_0$.
\item If in addition, we also assume that $p>\frac{d(3\beta-1)}{2(\beta-1)}$ and
$$
\kappa_1:=\sup_{t,\v}\|\Delta_x^{\frac{1}{3}}\sigma_t(\cdot,\v)\|_p+\|b\|_{\mH^{2/3,0}_p(T)}<\infty,
$$
then for any $f\in \mH^{2/3,0}_p(T)$, the unique solution $u$ also satisfies
\begin{align}\label{Pr1}
\|\nabla_x\nabla_\v u^\lambda\|_{\mL^p(T)}+\|\Delta^{\frac{1}{3}}_x\nabla^2_\v u^\lambda\|_{\mL^p(T)}\leq C\|f\|_{\mH^{2/3,0}_p(T)},
\end{align}
and for all $t\in[0,T]$,
\begin{align}\label{Pr2}
\|\Delta_x^{\frac{1}{3}}  u^\lambda_t\|_{\alpha,\beta;p}\leq C((T-t)\wedge\lambda^{-1})^{1-\frac{1}{p}-\frac{(3\alpha)\vee\beta}{2}}\|f\|_{\mH^{2/3,0}_p(t,T)},
\end{align}
where the constant $C$ only depends on $d,\delta,K,\alpha,\beta,p, T$ and $\kappa_1$.
\end{enumerate}
\et
\br
 In order to emphasize the dependence of the unique solution $u$ on $a,b$ and $T,\lambda, f$, 
we sometimes denote $u=\cR^{a,b}_{\lambda,T}(f)$.
\er
\subsection{Case $b=0$}
In this subsection we first consider the case of $b=0$ by using the freezing coefficient argument, and show the following basic existence and uniqueness result
for equation \eqref{PDE}.
\bt\label{Th33}
Let $p>1$ and $\alpha\in[0,\frac{2}{3})$, $\beta\in[0,2)$. Suppose {\bf (H$^{\delta,p}_K$)} holds. 
\begin{enumerate}[(i)]
\item For any $f\in \mL^p(T)$, there exists a unique solution $u=u^\lambda$ to \eqref{PDE} in the sense of Definition \ref{Def31} so that
\begin{align}\label{Pri}
\|u^\lambda\|_{\mH^{2/3,2}_p(T)}\leq C_1\|f\|_{\mL^p(T)}.
\end{align}
\item If $p>\frac{2}{(2-3\alpha)\wedge(2-\beta)}$, then for all $t\in[0,T]$,
\begin{align}\label{Pri2}
\|u^\lambda_t\|_{\alpha,\beta;p}\leq C_2((T-t)\wedge\lambda^{-1})^{1-\frac{1}{p}-\frac{(3\alpha)\vee\beta}{2}}\|f\|_{\mL^p(t,T)}.
\end{align}
\end{enumerate}
Here $C_1=C_1(d,\delta,K,p,T)$ and $C_2=C_2(d,\delta,K,\alpha,\beta,p,T)$ are increasing with respect to $T$.
\et
\begin{proof}
We show the apriori estimates \eqref{Pri} and \eqref{Pri2} by the freezing coefficient argument. 
The existence of a solution follows by the standard continuity argument. We divide the proof into five steps.

{\bf (a)} First of all, we assume that $u\in C([0,T]; H^{\infty,\infty}_p)$
satisfies \eqref{PDE} for Lebesgue almost all $t\in[0,T]$.
For given $p\geq 1$, let $\phi$ be a nonnegative symmetric smooth function on $\mR^{2d}$ with support in the unit ball and
$$
\int_{\mR^{2d}}|\phi(z)|^p\dif z=1.
$$
Let $\delta\in(0,1)$ be as in {\bf (H$^{\delta,p}_K$)} and set
$$
\phi_\delta(z):=\delta^{-2d/p}\phi(z/\delta),
$$
and for $z^o=(x^o,\v^o)$ and $t\in[0,T]$, define
$$
z^o_t:=(x^o-t\v^o,\v^o), \  \phi^{z^o_t}_\delta(z):=\phi_\delta(z^o_t-z).
$$
By definition, it is easy to see that
\begin{align}\label{LK1}
\int_{\mR^{2d}}|\phi^{z^o_t}_\delta(z)|^p\dif z^o=1,\ \ t\in[0,T],\ z\in\mR^{2d},
\end{align}
and for $j=1,2$,
\begin{align}\label{LK2}
\sup_{t\in[0,T]}\sup_{z\in\mR^{2d}}\int_{\mR^{2d}}|\nabla^j_\v\phi^{z^o_t}_\delta(z)|^p\dif z^o\leq C_\delta.
\end{align}
Define the freezing functions at point $z^o=(x^o,\v^o)$ as follows:
$$
a^{z^o}_t:=a_{t}(z^o_t),\ \ u^{z^o}_{\delta,t}(z):=u_t(z)\phi^{z^o_t}_\delta(z).
$$
By \eqref{PDE} and easy calculations, one sees that
\begin{align}\label{LK88}
\p_t u^{z^o}_\delta+\mathrm{tr}(a^{z^o}_t\cdot\nabla^2_{\v} u^{z^o}_\delta)+\v\cdot \nabla_x u^{z^o}_\delta-\lambda u^{z^o}_\delta=g^{z^o}_{\delta},
\end{align}
where for $z=(x,\v)$,
\begin{align*}
g^{z^o}_{\delta,t}(z)&:=\mathrm{tr}(a^{z^o}_t\cdot\nabla^2_{\v} u^{z^o}_{\delta,t})(z)-\mathrm{tr}(a_t\cdot\nabla^2_\v u_t)(z)\phi_\delta^{z^o_t}(z)\\
&\quad+(\v-\v^o)\cdot \nabla_x \phi_\delta^{z^o_t}(z) u_t(z)+f_t(z)\phi_\delta^{z^o_t}(z).
\end{align*}
We have the following {\it claim:} 
\begin{align}\label{LK3}
\begin{split}
&\left(\int_{\mR^{2d}}\|g^{z^o}_{\delta}\|_{\mL^p(t,T)}^p\dif z^o\right)^{1/p}
\leq \omega_a(\delta)\|\nabla^2_\v u\|_{\mL^p(t,T)}\\
&\quad+C_\delta\Big(\|u\|_{\mL^p(t,T)}+\|f\|_{\mL^p(t,T)}\Big),\ \ t\in[0,T].
\end{split}
\end{align}
{\it Proof of the claim:} Observe that
\begin{align*}
g^{z^o}_{\delta,t}(z)&=\mathrm{tr}((a^{z^o}_t-a_t)\cdot\nabla^2_\v u_t)(z)\phi_\delta^{z^o_t}(z)
+\mathrm{tr}(a^{z^o}_t\cdot(\nabla_\v u_t\otimes\nabla_\v\phi^{z^o_t}_\delta))(z)\\
&\quad+[\mathrm{tr}(a^{z^o}_t\cdot\nabla^2_\v\phi^{z^o_t}_\delta)(z)+(\v-\v^o)\cdot \nabla_x \phi_\delta^{z^o_t}(z)] u_t(z)
+f_t(z)\phi_\delta^{z^o_t}(z)\\
&=:I^\delta_1(t,z,z^o)+I^\delta_2(t,z,z^o)+I^\delta_3(t,z,z^o)+I^\delta_4(t,z,z^o).
\end{align*}
For $I^\delta_1(t,z,z^o)$, since the support of $\phi_\delta$ is in $B_\delta:=\{z\in\mR^{2d}: |z|\leq\delta\}$,
by  the definition of $\omega_a(\delta)$ and \eqref{LK1},  we have
\begin{align*}
\left(\int_{\mR^{2d}}\|I^\delta_1(\cdot,\cdot,z^o)\|_{\mL^p(t,T)}^p\dif z^o\right)^{1/p}
&\leq \omega_a(\delta)\left(\int_{\mR^{2d}}\|\nabla^2_\v u\cdot\phi^{z^o_\cdot}_\delta\|^p_{\mL^p(t,T)}\dif z^o\right)^{1/p}\\
&=\omega_a(\delta)\|\nabla^2_\v u\|_{\mL^p(t,T)}.
\end{align*}
For $I^\delta_2(t,z,z^o)$, by \eqref{LK2} we have
$$
\left(\int_{\mR^{2d}}\|I^\delta_2(\cdot,\cdot,z^o)\|_{\mL^p(t,T)}^p\dif z^o\right)^{1/p}\leq C_\delta\|\nabla_\v u\|_{\mL^p(t,T)}.
$$
For $I^\delta_3(t,z,z^o)$, we similarly have
$$
\left(\int_{\mR^{2d}}\|I^\delta_3(\cdot,\cdot,z^o)\|_{\mL^p(t,T)}^p\dif z^o\right)^{1/p}\leq C_\delta\|u\|_{\mL^p(t,T)}.
$$
For $I^\delta_4(t,z,z^o)$, by \eqref{LK1} we have
$$
\left(\int_{\mR^{2d}}\|I^\delta_4(\cdot,\cdot,z^o)\|_{\mL^p(t,T)}^p\dif z^o\right)^{1/p}=\|f\|_{\mL^p(t,T)}.
$$
Combining the above calculations, and by the interpolation inequality \eqref{Inter} and Young's inequality, we get the claim.

\

Now by \eqref{LK1}, we have
\begin{align}\label{ES6}
\begin{split}
&\|\nabla^2_\v u\|_{\mL^p(t,T)}
=\left(\int_{\mR^{2d}}\|(\nabla^2_\v u)\phi^{z^o_\cdot}_\delta\|_{\mL^p(t,T)}^p\dif z^o\right)^{1/p}\\
&\quad\leq \left(\int_{\mR^{2d}}\|\nabla^2_\v (u\phi^{z^o_\cdot}_\delta)-(\nabla^2_\v u)\phi^{z^o_\cdot}_\delta\|_{\mL^p(t,T)}^p\dif z^o\right)^{1/p}\\
&\qquad+\left(\int_{\mR^{2d}}\|\nabla^2_\v u^{z^o}_\delta\|_{\mL^p(t,T)}^p\dif z^o\right)^{1/p}=:I_1+I_2.
\end{split}
\end{align}
For $I_1$, by \eqref{LK2} and the interpolation inequality, we have
$$
I_1\leq C_\delta\Big(\|\nabla_\v u\|_{\mL^p(t,T)}+\|u\|_{\mL^p(t,T)}\Big)\leq \omega_a(\delta)\|\nabla^2_\v u\|_{\mL^p(t,T)}+C_{\delta}\|u\|_{\mL^p(t,T)}.
$$
For $I_2$, noticing that by \eqref{LK88} and Duhamel's formula,
\begin{align}\label{EQ1}
u^{z^o}_{\delta,t}(z)=\int^T_t\e^{\lambda(t-s)} P^{z^o}_{t,s}g^{z^o}_{\delta,s}(z)\dif s,
\end{align}
where $P^{z^o}_{t,s}$ is defined by \eqref{PP} in terms of $\sigma^{z^o}_t:=(a^{z^o}_t)^{1/2}$, 
by \eqref{Pri22} and \eqref{LK3}, we have 
\begin{align*}
\begin{split}
I_2\leq C_0&\left(\int_{\mR^{2d}}\|g^{z^o}_{\delta}\|_{\mL^p(t,T)}^p\dif z^o\right)^{1/p}
\leq C_0\omega_a(\delta)\|\nabla^2_\v u\|_{\mL^p(t,T)}\\
&\quad+C_\delta\Big(\|u\|_{\mL^p(t,T)}+C\|f\|_{\mL^p(t,T)}\Big).
\end{split}
\end{align*}
Substituting these two estimates into \eqref{ES6} and by $\omega_a(\delta)\leq\frac{1}{2(C_0+1)}$, 
we get
\begin{align}\label{LK9}
\|\nabla^2_\v u\|_{\mL^p(t,T)}\leq C\Big(\|u\|_{\mL^p(t,T)}+\|f\|_{\mL^p(t,T)}\Big).
\end{align}
Similarly, one can show that (see also step (c) below)
\begin{align}\label{LK8}
\|\Delta_x^{\frac{1}{3}}u\|_{\mL^p(t,T)}\leq C\Big(\|u\|_{\mL^p(t,T)}+\|f\|_{\mL^p(t,T)}\Big).
\end{align}

{\bf (b)}  By \eqref{EQ1} and the contraction of operator $P^{z^o}_{t,s}$ in $L^p(\mR^{2d})$, we have
\begin{align*}
\|u_t\|_p^p&\stackrel{\eqref{LK1}}{=}\int_{\mR^{2d}}\|u^{z^o}_{\delta,t}\|^p_p\dif z^o\leq\int_{\mR^{2d}}\left(\int^T_t \|g^{z^o}_{\delta,s}\|_{p}\dif s\right)^p\dif z^o\\
&\leq \left(\int^T_t\e^{p\lambda(t-s)/(p-1)}\dif s\right)^{p-1}\int_{\mR^{2d}}\int^T_t\|g^{z^o}_{\delta,s}\|^p_{p}\dif s\dif z^o\\
&\stackrel{\eqref{LK3}}{\preceq} ((T-t)\wedge\lambda^{-1})^{p-1}\Big(\|\nabla^2_\v u\|_{\mL^p(t,T)}^p+\|u\|_{\mL^p(t,T)}^p+\|f\|_{\mL^p(t,T)}^p\Big)\\
&\stackrel{\eqref{LK9}}{\preceq}((T-t)\wedge\lambda^{-1})^{p-1}\left(\int^T_t\|u_s\|_{p}^p\dif s+\|f\|_{\mL^p(t,T)}^p\right),
\end{align*}
which yields by Gronwall's inequality that
\begin{align}\label{ES77}
\|u_t\|^p_p\leq C((T-t)\wedge\lambda^{-1})^{p-1}\|f\|_{\mL^p(t,T)}^p.
\end{align}
Substituting it into \eqref{LK9} and \eqref{LK8}, we obtain \eqref{Pri}, and also by \eqref{LK3},
\begin{align}\label{LP2}
\int_{\mR^{2d}}\|g^{z^o}_{\delta}\|_{\mL^p(t,T)}^p\dif z^o\leq C\|f\|_{\mL^p(t,T)}^p,\ \ t\in[0,T].
\end{align}

{\bf (c)} Let $\alpha\in(0,\frac{2}{3})$.  By \eqref{EQ1}, \eqref{Gr} and H\"older's inequality, we have
\begin{align}
\|u^{z^o}_{\delta,t}\|_{\alpha,0; p}&\leq \int^T_t\e^{\lambda(t-s)} \|P^{z^o}_{t,s}g^{z^o}_{\delta,s}\|_{\alpha,0;p}\dif s
\preceq\int^T_t\e^{\lambda(t-s)}(s-t)^{-\frac{3\alpha}{2}}\|g^{z^o}_{\delta,s}\|_{p}\dif s\no\\
&\preceq \left(\int^T_t\e^{\lambda(t-s)} (s-t)^{-\frac{3p\alpha}{2(p-1)}}\dif s\right)^{1-\frac{1}{p}}\|g^{z^o}_{\delta}\|_{\mL^p(t,T)}\no\\
&\preceq( (T-t)\wedge\lambda^{-1})^{1-\frac{1}{p}-\frac{3\alpha}{2}}\|g^{z^o}_{\delta}\|_{\mL^p(t,T)}.\label{LP1}
\end{align}
By \eqref{LK1} again, we have
\begin{align}\label{LP4}
\begin{split}
\|\Delta_x^{\frac{\alpha}{2}}u_t\|^p_{p}&=\int_{\mR^{2d}}\|(\Delta_x^{\frac{\alpha}{2}} u_t)\phi^{z^o_t}_\delta\|^p_{p}\dif z^o 
\preceq\int_{\mR^{2d}}\|\Delta_x^{\frac{\alpha}{2}} u^{z^o}_{\delta,t}\|^p_{p}\dif z^o \\
&\qquad+\int_{\mR^{2d}}\|\Delta_x^{\frac{\alpha}{2}} (u_t\phi^{z^o_t}_\delta)-(\Delta_x^{\frac{\alpha}{2}} u_t)\phi^{z^o_t}_\delta\|^p_{p}\dif z^o.
\end{split}
\end{align}
By \eqref{LP1} and \eqref{LP2}, we have
\begin{align}
\begin{split}
\int_{\mR^{2d}}\|\Delta_x^{\frac{\alpha}{2}}u^{z^o}_{\delta,t}\|^p_{p}\dif z^o 
&\preceq((T-t)\wedge\lambda^{-1})^{p-1-\frac{3p\alpha}{2}}\int_{\mR^{2d}}\|g^{z^o}_{\delta}\|_{\mL^p(t,T)}^p\dif z^o\\
&\preceq((T-t)\wedge\lambda^{-1})^{p-1-\frac{3p\alpha}{2}}\|f\|^p_{\mL^p(t,T)}.
\end{split}\label{LP5}
\end{align}
On the other hand, noticing that by definition \eqref{Def1},
$$
\Delta_x^{\frac{\alpha}{2}} (u_t\phi^{z^o_t}_\delta)-(\Delta_x^{\frac{\alpha}{2}} u_t)\phi^{z^o_t}_\delta
=u_t\cdot\Delta_x^{\frac{\alpha}{2}}\phi^{z^o_t}_\delta+\int_{\mR^d}\delta_{(y,0)}u_t\cdot\delta_{(y,0)}\phi^{z^o_t}_\delta|y|^{-d-\alpha}\dif y,
$$
and
$$
\sup_z\int_{\mR^{2d}}|\Delta_x^{\frac{\alpha}{2}}\phi_\delta(z^o_t-z)|^p\dif z^o\leq C_\delta,
$$
$$
\sup_z\left(\int_{\mR^{2d}}|\delta_{(y,0)}\phi^{z^o_t}_\delta(z)|^p\dif z^o\right)^{1/p}\leq C_\delta(|y|\wedge 1),
$$
by Minkovskii's inequality, we have
\begin{align}\label{LP6}
\begin{split}
&\int_{\mR^{2d}}\|\Delta_x^{\frac{\alpha}{2}} (u_t\phi^{z^o_t}_\delta)-(\Delta_x^{\frac{\alpha}{2}} u_t)\phi^{z^o_t}_\delta\|_p^p\dif z^o
\preceq\int_{\mR^{2d}}\|u_t\cdot\Delta_x^{\frac{\alpha}{2}}\phi^{z^o_t}_\delta\|_p^p\dif z^o\\
&\qquad+\int_{\mR^{2d}}\left(\int_{\mR^d}|\delta_{(y,0)}u_t(z)|\cdot\left(\int_{\mR^{2d}}|\delta_{(y,0)}\phi^{z^o_t}_\delta(z)|^p\dif z^o\right)^{1/p}|y|^{-d-\alpha}\dif y\right)^p\dif z\\
&\qquad\preceq \|u_t\|^p_p+\int_{\mR^{2d}}\left(\int_{\mR^d}|\delta_{(y,0)}u_t(z)|\cdot(|y|\wedge 1)|y|^{-d-\alpha}\dif y\right)^p\dif z\\
&\qquad\preceq \|u_t\|^p_p+\left(\int_{\mR^d}\|\delta_{(y,0)}u_t\|_p\cdot(|y|\wedge 1)|y|^{-d-\alpha}\dif y\right)^p\preceq\|u_t\|^p_p.
\end{split}
\end{align}
Combining \eqref{LP4}, \eqref{LP5} and \eqref{LP6} with \eqref{ES77}, we arrive at
$$
\|u_t\|^p_{\alpha,0; p}\preceq((T-t)\wedge\lambda^{-1})^{p-1-\frac{3p\alpha}{2}}\|f\|^p_{\mL^p(t,T)}.
$$ 
Similarly, for any $\beta\in(0,2)$, one can show that
$$
\|u_t\|^p_{0,\beta; p}\preceq((T-t)\wedge\lambda^{-1})^{p-1-\frac{p\beta}{2}}\|f\|^p_{\mL^p(t,T)}.
$$ 
Combining the above two estimates, we obtain \eqref{Pri2}.

\

{\bf (d)} Below we assume that $u\in\mH^{0,2}_p(T)$ is a solution of \eqref{PDE} in the sense of Definition \ref{Def31}.
Let $\varrho_\eps$ be defined by \eqref{Rho} and set
$$
u_\eps:=u*\varrho_\eps,\ \ f_\eps:=f*\varrho_\eps. 
$$
Taking $\varphi=\varrho_\eps(z-\cdot)$ in \eqref{Sol}, we obtain
$$
\p_t u_\eps+\sL^a_t u_\eps-\lambda u_\eps+[\varrho_\eps,\sL^a_t]u_t+f_\eps=0,\ \ u_{\eps,T}=0,
$$
where $\sL^a_t:=\sL^{a,0}_t$ and $[\varrho_\eps,\sL^a_t]$ is defined by \eqref{Def0}. By what we have proved, it holds that
$$
\|u_\eps\|_{\mH^{2/3,2}_p(T)}\leq C_1(\|f_\eps\|_{\mL^p(T)}+\|[\varrho_\eps,\sL^{a}]\|_{\mL^p(T)}).
$$
Since $\nabla^2_\v u\in\mL^p(T)$, by the property of convolutions and (ii) of Lemma \ref{Le24}, we get \eqref{Pri} by taking limits.   Similarly, we also have \eqref{Pri2}.

\

{\bf (e)} Finally, we use the standard continuity argument to show the existence of a solution (see \cite{Kr0}). Consider the following parametrized equation:
\begin{align}\label{PPE}
\p_t u+\sL^{a_\tau}_t u-\lambda u+f=0, \ \ u_T=0,
\end{align}
where $\tau\in[0,1]$ and $a_\tau:=K(1-\tau)\mI+\tau a$. Since $K^{-1}\cdot\mI\leq a\leq K\cdot\mI$, we obviously have
$$
K^{-1}\cdot\mI\leq a_\tau\leq K\cdot \mI,\ \ \omega_{a_\tau}(\delta)=\omega_a(\delta).
$$ 
Hence the apriori estimate \eqref{Pri} holds for \eqref{PPE} with constant $C_1$ independent of $\tau\in[0,1]$.
Suppose that \eqref{PPE} is solvable for some $\tau_0\in[0,1)$. We want to show that \eqref{PPE} is also solvable for any $\tau\in[\tau_0,\tau_0+\frac{1}{4KC_1})$, 
where $C_1$ is the constant in \eqref{Pri}.
Let $u^0=0$ and for $n\in\mN$, define $u^n$ recursively by
$$
\p_t u^{n}+\sL^{a_{\tau_0}}_t u^{n}-\lambda u^n+\tr((a_\tau-a_{\tau_0})\cdot\nabla^2_\v u^{n-1})+f=0, \ \ u^{n}_T=0.
$$
By the apriori estimate \eqref{Pri}, we have
\begin{align*}
\|u^n\|_{\mH^{2/3,2}_p(T)}&\leq C_1\|\tr((a_\tau-a_{\tau_0})\cdot\nabla^2_\v u^{n-1})+f\|_{\mL^p(T)}\\
&\leq 2KC_1(\tau-\tau_0)\|u^{n-1}\|_{\mH^{2/3,2}_p(T)}+C_1\|f\|_{\mL^p(T)}\\
&\leq \tfrac{1}{2}\|u^{n-1}\|_{\mH^{2/3,2}_p(T)}+C_1\|f\|_{\mL^p(T)},
\end{align*}
and similarly,
$$
\|u^n-u^m\|_{\mH^{2/3,2}_p(T)}\leq \tfrac{1}{2}\|u^{n-1}-u^{m-1}\|_{\mH^{2/3,2}_p(T)},
$$
which imply that $u^n$ is a Cauchy sequence in $\mH^{2/3,2}_p(T)$. It is easy to see that the limit $u$ of $u^n$ satisfies
\eqref{PPE}. Since \eqref{PPE} is solvable for $\tau=0$ by (ii) of Theorem \ref{Th25}, 
by repeatedly using what we have proved finitely many times, we get the solvability of \eqref{PPE} for $\tau=1$.
\end{proof}
Next we show further regularity of the solution under extra assumption.
\bt
Suppose that for some $p>3d/2$,  {\bf (H$^{\delta,p}_K$)} holds and
$$
\kappa_1:=\sup_{(t,\v)\in[0,T]\times\mR^d}\|\Delta_x^{\frac{1}{3}}\sigma_t(\cdot,\v)\|_p<\infty.
$$
\begin{enumerate}[(i)]
\item For any $f\in \mH^{2/3,0}_p(T)$, the unique solution $u$ of PDE \eqref{PDE} also satisfies
\begin{align}\label{E35}
\|\nabla_x\nabla_\v u\|_{\mL^p(T)}+\|\Delta^{\frac{1}{3}}_x\nabla^2_\v u\|_{\mL^p(T)}\leq C_3\|f\|_{\mH^{2/3,0}_p(T)}.
\end{align}
\item Let $\alpha\in[0,\frac{2}{3})$ and $\beta\in[0,2)$. If $p>\frac{2}{(2-3\alpha)(2-\beta)}\vee\frac{3d}{2}$, then
\begin{align}\label{E36}
\|\Delta_x^{\frac{1}{3}} u_t\|_{\alpha,\beta;p}\leq C_4((T-t)\wedge\lambda^{-1})^{1-\frac{1}{p}-\frac{(3\alpha)\vee\beta}{2}}\|f\|_{\mH^{2/3,0}_p(T)}.
\end{align}
\end{enumerate}
Here $C_3=C_3(d,\delta,K,\kappa_1, p,T)$ and $C_4=C_4(d,\delta,K,\kappa_1, \alpha,\beta,p,T)$ are increasing with respect to $T$.
\et
\begin{proof}
As in the proof of Theorem \ref{Th33}, it suffices to show the apriori estimates \eqref{E35} and \eqref{E36}.
Notice that using \eqref{E28} with $\alpha=\frac{2}{3}$, $\beta=2$ and by \eqref{Pri},
\begin{align}\label{ES00}
\|\Delta_x^{\frac{1}{6}}\nabla_\v u\|_{\mL^p(t,T)}\preceq\|u\|_{\mH^{2/3,2}_p(t,T)}\preceq\|f\|_{\mL^p(t,T)}.
\end{align}
Assume $u\in C([0,T];H^{\infty,\infty}_p)$ and let $w_t(x,\v):=\Delta_x^{\frac{1}{3}}u_t(x,\v)$. By \eqref{PDE} we have
$$
\p_t w+\sL^a_t w+[\Delta_x^{\frac{1}{3}},\sL^a_t]u+\Delta_x^{\frac{1}{3}}f=0, w_T=0.
$$
By definition \eqref{Def1} and the assumptions, it is easy to see that
$$
\sup_{(t,\v)\in[0,T]\times\mR^d}\|\Delta_x^{\frac{1}{3}}a_t(\cdot,\v)\|_p<\infty.
$$
Hence, by \eqref{ES00}, \eqref{Pri} and \eqref{E210}, we have for any $\eps>0$,
\begin{align*}
&\|\nabla_x\nabla_\v u\|_{\mL^p(T)}+\|\Delta_x^{\frac{2}{3}}u\|_{\mL^p(T)}+\|\Delta_x^{\frac{1}{3}}\nabla^2_\v u\|_{\mL^p(T)}\\
&\quad\stackrel{\eqref{Riesz}}{\preceq}\|\Delta_x^{\frac{1}{6}}\nabla_\v w\|_{\mL^p(T)}+\|\Delta_x^{\frac{1}{3}} w\|_{\mL^p(T)}+\|\nabla^2_\v w\|_{\mL^p(T)}\\
&\quad\preceq \|[\Delta_x^{\frac{1}{3}},\sL^a_t]u+\Delta_x^{\frac{1}{3}}f\|_{\mL^p(T)}\\
&\quad\leq\|[\Delta_x^{\frac{1}{3}},\sL^a_t]u\|_{\mL^p(T)}+\|f\|_{\mH^{2/3,0}_p(T)}\\
&\quad\preceq\eps\|\Delta_x^{\frac{1}{3}}\nabla^2_\v u\|_{\mL^p(T)}+C_\eps\big(\|\nabla^2_\v u\|_{\mL^p(T)}+\|f\|_{\mH^{2/3,0}_p(T)}\big)\\
&\quad\preceq\eps\|\Delta_x^{\frac{1}{3}}\nabla^2_\v u\|_{\mL^p(T)}+C_\eps\|f\|_{\mH^{2/3,0}_p(T)},
\end{align*}
which implies the desired estimate \eqref{E35} by letting $\eps$ be small enough. As for \eqref{E36}, it follows by applying \eqref{Pri2} to $w$ and using the above estimate.
\end{proof}
\subsection{Proof of Theorem \ref{Th32}}
By a standard fixed point argument or Picard's iteration, it suffices to prove the apriori estimate \eqref{Pri0}. 
Let $\alpha\in(0,\frac{2}{3})$, $\beta\in(1,2)$ and $p>\frac{2}{(2-3\alpha)\wedge(2-\beta)}$ be not equal to $\frac{d(\alpha+\beta)}{\alpha(\beta-1)}$.

(i) Let $q\in[p\vee \frac{d(\alpha+\beta)}{\alpha(\beta-1)},\infty]$ and $r\in[p,\infty]$ with $\frac{1}{p}=\frac{1}{q}+\frac{1}{r}$. 
Notice that by H\"older's inequality and \eqref{E28}, and using $(\alpha(\beta-1)/\beta, \beta-1)$ in place of $(\alpha,\beta)$ in \eqref{E29},
$$
\|b\cdot\nabla_\v u\|^p_{\mL^p(t,T)}\leq \int^T_t\|b_s\|^p_{q}\|\nabla_\v u_s\|^p_{r}\dif s\preceq \int^T_t\|b_s\|^p_{q}\|u_s\|^p_{\alpha,\beta; p}\dif s.
$$
Thus, by \eqref{Pri2} we have
\begin{align*}
\|u_t\|_{\alpha,\beta;p}^p&\preceq ((T-t)\wedge\lambda^{-1})^{p-1-\frac{p((3\alpha)\vee\beta)}{2}}\|b\cdot\nabla_\v u+f\|^p_{\mL^p(t,T)}\\
&\preceq((T-t)\wedge\lambda^{-1})^{p-1-\frac{p((3\alpha)\vee\beta)}{2}}\|f\|^p_{\mL^p(t,T)}+\int^T_t\|b_s\|^p_{q}\|u_s\|^p_{\alpha,\beta;p}\dif s,
\end{align*}
which implies \eqref{Pri00} by Gronwall's inequality. 

On the other hand, by \eqref{Pri} we have
$$
\|u\|_{\mH^{2/3,0}_p(T)}^p\preceq \|b\cdot\nabla_\v u+f\|^p_{\mL^p(T)}\preceq
\int^T_0\|b_s\|^p_{q}\|u_s\|^p_{\alpha,\beta;p}\dif s+\|f\|^p_{\mL^p(T)},
$$
which in turn implies \eqref{Pri0} by \eqref{Pri00}.

(ii) Let $p>\frac{d(3\beta-1)}{2(\beta-1)}$. By Lemma \ref{Le23} with $\alpha=\frac{2}{3}$, we have
\begin{align*}
\|b\cdot\nabla_\v u\|^p_{\mH^{2/3,0}_p(t,T)}
&\preceq \int^T_t\|b_s\|^p_{2/3,0;p}\Big(\|\Delta_x^{\frac{1}{3}}\nabla_\v u_s\|_{0,\beta-1; p}^p+\|\nabla_\v u_s\|_{0,\beta-1;p}^p\Big)\dif s\\
&\preceq \int^T_t\|b_s\|^p_{2/3,0;p}\Big(\|\Delta_x^{\frac{1}{3}}u_s\|_{0,\beta; p}^p+\|u_s\|_{0,\beta;p}^p\Big)\dif s.
\end{align*}
Thus, by \eqref{E36} and \eqref{Pri00}, we have
\begin{align*}
\|\Delta_x^{\frac{1}{3}} u_t\|^p_{\alpha,\beta;p}&\preceq ((T-t)\wedge\lambda^{-1})^{p-1-\frac{p((3\alpha)\vee\beta)}{2}}\|b\cdot\nabla_\v u+f\|^p_{\mH^{2/3,0}_p(t,T)}\\
&\preceq ((T-t)\wedge\lambda^{-1})^{p-1-\frac{p((3\alpha)\vee\beta)}{2}}\|f\|^p_{\mH^{2/3,0}_p(t,T)}\\
&\quad+\int^T_t\|b_s\|^p_{2/3,0;p}\|\Delta_x^{\frac{1}{3}}u_s\|_{0,\beta; p}^p\dif s,
\end{align*}
which yields  \eqref{Pr2} by Gronwall's inequality.

Moreover, by \eqref{E35} we have
\begin{align*}
\|\nabla_x\nabla_\v u\|_{\mL^p(T)}+\|\Delta_x^{\frac{1}{3}}\nabla^2_\v u\|_{\mL^p(T)}\leq 
C\|b\cdot\nabla_\v u+f\|_{\mH^{2/3,0}_p(T)}\leq C\|f\|_{\mH^{2/3,0}_p(T)},
\end{align*}
which gives \eqref{Pr1}. The proof is complete.

\section{Well-posedness of martingale problem}
Let $\Omega=C(\mR_+;\mR^{2d})$ be the space of all continuous functions from $\mR_+$ to $\mR^{2d}$, which is endowed 
with the locally uniform convergence topology. Let 
$$
Z_t(\omega):=\omega_t
$$
be the coordinate process on $\Omega$. For $t\geq 0$, let
$$
\sF_t:=\sigma\Big\{Z_s: s\in[0,t]\Big\}, \ \ \sF:=\vee_{t\geq 0}\sF_t.
$$
We first recall the following notions of martingale solution and weak solution (see \cite{St-Va}).
\bd
Let $\sigma:\mR_+\times\mR^{2d}\to\mR^d\otimes\mR^d$ and $b:\mR_+\times\mR^{2d}\to\mR^d$ be Borel measurable functions and
$\sL^{a,b}_t$ be defined by \eqref{LL} with $a=\frac{1}{2}\sigma\sigma^*$.
\begin{enumerate}[(i)]
\item (Martingale solution) For given $(r,z)\in\mR_+\times\mR^{2d}$,  a probability measure $\mP=\mP_{r,z}$ 
on $(\Omega,\sF)$ is said a solution to the martingale problem for $\sL^{a,b}_t$ starting from $(r,z)$ if 
$$
\mP\left(Z_t=z,t\in[0,r]\ \mathrm{ and } \int^t_r(\|a_s(Z_s)\|+|b_s(Z_s)|)\dif s<\infty,t\geq r\right)=1,
$$
and for all $\varphi\in C^\infty_c(\mR^d)$,
\begin{align}\label{MM4}
[r,\infty)\ni t\mapsto \varphi(Z_t)-\int^t_r\sL^{a,b}_s \varphi(Z_s)\dif s=:M^{r,\varphi}_t
\end{align}
is an $\sF_t$-martingale with respect to $\mP$ after time $r$. 
We denote by $\sP^{\sigma,b}_{r,z}$ the set of all martingale solutions associated with $(\sigma,b)$ and starting from $(r,z)$.
\item (Well-posedness) One says that the martingale problem for $\sL^{a,b}_t$ is well-posed if for each $(r,z)\in\mR_+\times\mR^{2d}$, there is exactly
one solution $\mP_{r,z}$ to the martingale problem  for $\sL^{a,b}_t$ starting from $(r,z)$.
\item (Weak solution) A triple $((\tilde Z,\tilde W); (\tilde \Omega,\tilde \sF,\tilde \mP);(\tilde \sF_t)_{t\geq 0})$ is called a weak solution of SDE \eqref{SDE0} 
with starting point $(r,z)\in\mR_+\times\mR^{2d}$ if $(\tilde \sF_t)_{t\geq 0}$ satisfies the usual conditions,
$\tilde W$ is a $d$-dimensional $\tilde \sF_t$-Brownian motion, and $\tilde Z=(\tilde X,\dot{\tilde  X})$ is an $\tilde \sF_t$-adapted $\mR^{2d}$-valued process satisfying
that 
$$
\int^t_r(\|a_s(\tilde Z_s)\|+|b_s(\tilde Z_s)|)\dif s<\infty, \ t\geq r,\ \  \tilde\mP-a.s.,
$$
and
$$
\tilde Z_t=z+\int^t_r(\dot{\tilde  X}_s,b_s(\tilde Z_s))\dif s+\int^t_r(0,\sigma_s(\tilde Z_s)\dif\tilde  W_s).
$$
\end{enumerate}
\ed
\br
It is well known that the martingale solutions and weak solutions are equivalent, for example, see \cite[p.318, Proposition 4.11]{Ka-Sh}.
\er
\subsection{Krylov's type estimate}
In this  subsection we first show the following important  estimate of Krylov's type for  weak solutions.
\bt\label{Th41}
Suppose that $\sigma$ satisfies {\bf (UE)} and
\begin{align}\label{EJ1}
\lim_{|z-z'|\to 0}\sup_{t\in[0,T]}\|\sigma_t(z)-\sigma_t(z')\|=0,\  T>0,
\end{align}
and $b\in \cap_T\mL^q(T)$ for some $q>(2(2d+1),\infty]$. 
Then for any $p>2d+1$ and $T>0$, there is a constant $C>0$ depending on $T,\sigma,p,q,d, \|b\|_{\mL^q(T)}$ 
such that for any $(r,z)\in[0,T)\times\mR^{2d}$, any weak solution 
$((\tilde Z,\tilde W); (\tilde \Omega,\tilde \sF,\tilde \mP);(\tilde \sF_t)_{t\geq 0})$ of SDE \eqref{SDE0} 
with starting point $(r,z)$, and $r\leq t_0<t_1\leq T$ and $f\in\mL^p(T)$,
\begin{align}\label{BB9}
\tilde\mE\left(\int^{t_1}_{t_0} f_s(\tilde Z_s)\dif s\Big|\tilde \sF_{t_0}\right)
\leq C(t_1-t_0)^{\frac{1}{2d+1}-\frac{1}{p}}\|f\|_{\mL^p(t_0,t_1)},
\end{align}
where $C>0$ is increasing in $T,\omega_\sigma(\delta)$ and $\|b\|_{\mL^q(T)}$.
\et
\begin{proof}
Without loss of generality, we assume $(r,z)=(0,0)$ and drop the tilde in the definition of weak solutions  for simplicity. 
We divide the proof into four steps.

{\bf (a)} Let $a_t(z):=\frac{1}{2}(\sigma_t\sigma^*_t)(z)$ and $p\in(2(2d+1),q]$.
For any $0\leq t_0<t_1\leq T$, $\lambda\geq 1$ and $f\in\mL^p(t_0,t_1)$,
let $u=\cR^{a,b}_{\lambda,t_1}(f)$ be the solution of PDE \eqref{PDE} with terminal time $T=t_1$.
By \eqref{E29} and \eqref{Pri00} with $\alpha=\frac{4d}{3(2d+1)}$ and $\beta=\frac{4d}{2d+1}$,
\begin{align}\label{44}
\begin{split}
&\|\cR^{a,b}_{\lambda,t_1}(f)\|_{\mL^\infty(t_0,t_1)}
\preceq \big\|\cR^{a,b}_{\lambda,t_1}(f)\big\|_{L^\infty([t_0,t_1];H^{\alpha,\beta}_p)}\\
&\quad\preceq ((t_1-t_0)\wedge\lambda^{-1})^{\frac{1}{2d+1}-\frac{1}{p}}\|f\|_{\mL^p(t_0,t_1)},
\end{split}
\end{align}
and by \eqref{E28}, \eqref{E29},  \eqref{Pri00} with $\alpha=\frac{4d+1}{3(2d+1)}$ and $\beta=\frac{4d+1}{2d+1}$,
\begin{align}\label{404}
\begin{split}
&\|\cR^{a,0}_{\lambda,t_1}(f)\|_{\mL^\infty(t_0,t_1)}\preceq 
\big\|\cR^{a,0}_{\lambda,t_1}(f)\big\|_{L^\infty([t_0,t_1];H^{\alpha,\beta}_p)}\\
&\quad\preceq ((t_1-t_0)\wedge\lambda^{-1})^{\frac{1}{2(2d+1)}-\frac{1}{p}}\|f\|_{\mL^p(t_0,t_1)},
\end{split}
\end{align}
where the $C$ in the above $\preceq$ only depend on $d,p,q,K,T,\omega_a(\delta)$ and $\|b\|_{\mL^q(T)}$.

Now, for any $R>0$, define a stopping time
\begin{align}\label{TR}
\tau_R:=\Bigg\{t\geq 0: \int^t_0 |b_s(Z_s)|\dif s>R\Bigg\}.
\end{align}
Let $\varrho_\eps$ be as in \eqref{Rho}. We introduce a $d\times d$ matrix-valued function, which is crucial for us below. For 
$t\geq 0$ and $z\in\mR^{2d}$, let
$$
a^{\eps,R}_{t}(z):=
\left\{
\begin{array}{ll}
\tfrac{\mE\big(a_t(Z_t)\varrho_\eps(Z_t-z)1_{t<\tau_R}\big)}{\mE\big(\varrho_\eps(Z_t-z)1_{t<\tau_R}\big)},\ 
& \mathrm{ if }\ \ \mE\big(\varrho_\eps(Z_t-z)1_{t<\tau_R}\big)\not=0,\\
a_t(z),\ & \mathrm{ if }\ \ \mE\big(\varrho_\eps(Z_t-z)1_{t<\tau_R}\big)=0.
\end{array}
\right.
$$
By {\bf (UE)} and the definitions, we have
$$
K^{-2}\cdot\mI\leq a^{\eps,R}_{t}(z)\leq K^2\cdot\mI,
$$
and for all $|z-z'|\leq\eps$,
\begin{align*}
\|a^{\eps,R}_{t}(z)-a^{\eps,R}_{t}(z')\|&\leq\|a_{t}(z)-a_{t}(z')\|+\left\|\tfrac{\mE\big((a_t(Z_t)-a_t(z))\varrho_\eps(Z_t-z)1_{t<\tau_R}\big)}
{\mE\big(\varrho_\eps(Z_t-z)1_{t<\tau_R}\big)}\right\|\\
&+\left\|\tfrac{\mE\big((a_t(Z_t)-a_t(z'))\varrho_\eps(Z_t-z')1_{t<\tau_R}\big)}{\mE\big(\varrho_\eps(Z_t-z)1_{t<\tau_R}\big)}\right\|\leq 3\omega_a(\eps).
\end{align*}
Let $C_0$ be the same as in \eqref{Pri22}. By \eqref{EJ1}, one may choose $\delta_0$ small enough such that for all $\eps\in(0,\delta_0)$ and $R>0$,
\begin{align}\label{App}
\omega_{a^{\eps,R}}(\eps)\leq 3\omega_a(\delta_0)\leq\tfrac{1}{2(C_0+1)}.
\end{align}

{\bf (b)} In this step we show that for any $p>2(2d+1)$  and $f\in\mL^p(T)$,
\begin{align}\label{BB90}
\mE\left(\int^{T}_0 f_s(Z_s)\dif s\right)\leq C\|f\|_{\mL^p(T)}.
\end{align}
By a standard density argument, we may assume $f\in C_c([0,T]\times\mR^{2d})$. Let 
$$
u^{\eps,R}:=\cR^{a^{\eps,R},0}_{\lambda, T}(f),\ \ u^{\eps,R}_{\eps}:=u^{\eps,R}*\varrho_\eps,\ \ f_\eps:=f*\varrho_\eps.
$$
By It\^o's formula, we have
\begin{align}\label{EY1}
\begin{split}
&\mE\Big(u^{\eps,R}_{\eps,T\wedge\tau_R}(Z_{T\wedge\tau_R})-u_{\eps,0}(Z_0)\Big)
=\mE\left(\int^{T\wedge\tau_R}_0\!\!\! (\p_s u^{\eps,R}_{\eps,s}+\sL^{a,b}_su^{\eps,R}_{\eps,s})(Z_s)\dif s\right).
\end{split}
\end{align}
Noticing that by Definition \ref{Def31},
$$
\p_su^{\eps,R}_{\eps, s}+(\sL^{a^{\eps,R},0}_su^{\eps,R}_s)*\varrho_\eps-\lambda u^{\eps,R}_{\eps, s}+f_\eps=0,
$$
and by the definitions of $u^{\eps,R}_\eps$ and $a^{\eps,R}$,
\begin{align*}
&\mE\big(\tr(a_s\cdot\nabla^2_\v u^{\eps,R}_{\eps,t})(Z_s)1_{s<\tau_R}\big)\\
&=\int_{\mR^{2d}}\mE\big(\tr(a_s(Z_s)\cdot\nabla^2_\v u^{\eps,R}_s(z))\varrho_\eps(Z_s-z)1_{s<\tau_R}\big)\dif z\\
&=\int_{\mR^{2d}}\tr(a^{\eps,R}_s(z)\cdot\nabla^2_\v u^{\eps,R}_s(z))\mE\big(\varrho_\eps(Z_s-z)1_{s<\tau_R}\big)\dif z\\
&=\mE\Big((\tr(a^{\eps,R}_s\cdot\nabla^2_\v u^{\eps,R}_s)*\varrho_\eps)(Z_s)1_{s<\tau_R}\Big),
\end{align*}
by easy calculations, one sees that
\begin{align*}
&\mE\Big((\p_s u^{\eps,R}_{\eps,t}+\sL^{a,b}_su^{\eps,R}_{\eps,t})(Z_s)1_{s<\tau_R}\Big)\\
&=\int_{\mR^{2d}}\!\!\!\mE\Big((\dot X_s-\v)\cdot\nabla_x\varrho_\eps(Z_s-z)1_{s<\tau_R}\Big)u^{\eps,R}_s(z)\dif z\\
&\quad+\mE\Big((b_s\cdot\nabla_\v u^{\eps,R}_{\eps,t}+\lambda u^{\eps,R}_s-f_\eps)(Z_s)1_{s<\tau_R}\Big)\\
&\leq\|u^{\eps,R}_s\|_\infty\left(\lambda+\int_{\mR^{2d}}|\v|~|\nabla_x\varrho_\eps|(x,\v)\dif x\dif\v\right)\\
&\quad+\|\nabla_\v u^{\eps,R}_s\|_\infty\mE\big(|b_s|(Z_s)1_{s<\tau_R}\big)
-\mE\big( f_\eps(Z_s)1_{s<\tau_R}\big).
\end{align*}
Substituting this into \eqref{EY1}, we obtain
\begin{align*}
&\mE\left(\int^{T\wedge\tau_R}_0f_\eps(Z_s)\dif s\right)\leq 
(\|\nabla_x\varrho\|_1+2+\lambda)\|u^{\eps,R}\|_{\mL^\infty(T)}\\
&\qquad+\|\nabla_\v u^{\eps,R}\|_{\mL^\infty(T)}\mE\left(\int^{T\wedge\tau_R}_0|b_s|(Z_s)\dif s\right).
\end{align*}
By \eqref{App}, \eqref{44} and \eqref{404} with $b\equiv0$, there is a $C=C(d,p,K,T,\omega_a(\delta_0))>0$  such that for all $\eps\in(0,\delta_0)$ and $\lambda>1$,
\begin{align*}
&\mE\left(\int^{T\wedge\tau_R}_0f_{\eps,s}(Z_s)\dif s\right)\leq 
C(1+\lambda)T^{\frac{1}{2d+1}-\frac{1}{p}}\|f\|_{\mL^p(T)}\\
&\qquad+C\lambda^{\frac{1}{p}-\frac{1}{2(2d+1)}}\|f\|_{\mL^p(T)}\mE\left(\int^{T\wedge\tau_R}_0|b_s|(Z_s)\dif s\right),
\end{align*}
which  implies, by letting $\eps\to 0$ and $\lambda$ large enough, that for any $\delta>0$, there is a $C_\delta>0$ such that for all $f\in\mL^p(T)$,
\begin{align}\label{BB090}
\mE\left(\int^{T\wedge\tau_R}_0 f_s(Z_s)\dif s\right)\leq \left(C_\delta+\delta\mE\left(\int^{T\wedge\tau_R}_0|b_s|(Z_s)\dif s\right)\right)\|f\|_{\mL^p(T)}.
\end{align}
In particular, choosing $f_s=|b_s|$ and $\delta\leq1/(2\|b\|_{\mL^q(T)})$, we get
$$
\mE\left(\int^{T\wedge\tau_R}_0|b_s|(Z_s)\dif s\right)\leq C\|b\|_{\mL^q(T)}.
$$
Substituting this into \eqref{BB090} and letting $R\to\infty$, we get \eqref{BB90}.

{\bf (c)} In this step we show that \eqref{BB9} holds for $p=q>2(2d+1)$. Let $0\leq t_0<t_1\leq T$ and $f\in \mL^q(t_0,t_1)$, and write
$$
u:=-\cR^{a,b}_{0,t_1}(f)\in\mH^{0,2}_q(t_1),\ \ u_\eps:=u*\varrho_\eps.
$$
Noticing that by definitions,
\begin{align*}
\p_t u_{\eps,t}+\sL^{a,b}_tu_{\eps,t}=f_{\eps,t}+[\varrho_\eps,\sL^{a,b}_t]u_t,\ \ u_{\eps,t_1}=0,
\end{align*} 
by It\^o's formula, we have
\begin{align}
\mE^{\sF_{t_0}}\Big(u_{\eps,t_1}(Z_{t_1})-u_{\eps,t_0}(Z_{t_0})\Big)
=\mE^{\sF_{t_0}}\!\left(\int^{t_1}_{t_0} \Big(f_{\eps,s}+[\varrho_\eps,\sL^{a,b}_s]u_s\Big)(Z_s)\dif s\right),\label{BB2}
\end{align}
where $\mE^{\sF_{t_0}}(\cdot)=\mE(\cdot|^{\sF_{t_0}})$.
Since by \eqref{BB90}, \eqref{404} and Lemma \ref{Le24},
$$
\lim_{\eps\to 0}\mE\left(\int^{t_1}_{t_0}|[\varrho_\eps,\sL^{a,b}_s]u_s|(Z_s)\dif s\right)\preceq
\lim_{\eps\to 0}\|[\varrho_\eps,\sL^{a,b}_\cdot]u\|_{\mL^q(t_0,t_1)}=0
$$
and
$$
\lim_{\eps\to 0}\mE\left(\int^{t_1}_{t_0}|f_{\eps,s}-f_s|(Z_s)\dif s\right)\preceq
\lim_{\eps\to 0}\|f_\eps-f\|_{\mL^q(t_0,t_1)}=0,
$$
taking limits $\eps\to 0$ for both sides of (\ref{BB2}) and by \eqref{44}, we obtain
\begin{align}\label{KJ0}
\mE^{\sF_{t_0}}\left(\int^{t_1}_{t_0} f_s(Z_s)\dif s\right)
\leq 2\|u\|_{\mL^\infty(t_0,t_1)}\leq C(t_1-t_0)^{\frac{1}{2d+1}-\frac{1}{q}}\|f\|_{\mL^q(t_0,t_1)}.
\end{align}

{\bf (d)} Let $0\leq t_0<t_1\leq T$. By \eqref{KJ0} with $f_s=|b_s|$ and Corollary \ref{Cor1} below, we have for any $\lambda>0$,
\begin{align}\label{KJ1}
\mE^{\sF_{t_0}}\exp\left(\lambda\int^{t_1}_{t_0}|b_s(Z_s)|\dif s\right)\leq C(\lambda,\|b\|_{\mL^q(T)})<\infty.
\end{align}
Define
$$
\cE_{t_0,t_1}:=\exp\left\{\int^{t_1}_{t_0}(\sigma^{-1}_sb_s)(Z_s)\dif W_s-\frac{1}{2}\int^{t_1}_{t_0}|(\sigma^{-1}_sb_s)(Z_s)|^2\dif s\right\}.
$$
By Novikov's criterion, $\mE(\cE_{t_0,t_1})=1$, and for any $\gamma\in\mR$, by \eqref{KJ1} and H\"older's inequality,
\begin{align}\label{KJ2}
\mE\Big(\cE_{t_0,t_1}^\gamma|\sF_{t_0}\Big)<\infty.
\end{align}
Define a new probability $\mQ_{t_0, t_1}:=\cE_{t_0,t_1}\mP$. By Girsanov's theorem, $\bar W_\cdot:=W_\cdot+\int^\cdot_{t_0}b_s(Z_s)\dif s$ 
is still a Brownian motion under $\mQ_{t_0,t_1}$.
Moreover, $Z_t$ satisfies
$$
Z_t=Z_{t_0}+\int^t_{t_0}(\dot{X}_s,0)\dif s+\int^t_{t_0}(0,\sigma_s(Z_s)\dif\bar  W_s).
$$
By the same argument as used in {\bf (c)} with $b\equiv 0$, since in this case, we only need to control $\|u\|_{\mL^\infty(t_0,t_1)}$ and \eqref{44} holds for any $p>2d+1$, 
we obtain that there
is a constant $C>0$ such that for all $p>2d+1$ and $f\in\mL^p(T)$, 
$$
\mE^{\mQ_{t_0,t_1}}\left(\int^{t_1}_{t_0} f_s(Z_s)\dif s\Big|\sF_{t_0}\right)
\leq C(t_1-t_0)^{\frac{1}{2d+1}-\frac{1}{p}}\|f\|_{\mL^p(t_0,t_1)}.
$$
Finally, by \eqref{KJ2} and H\"older's inequality, we obtain \eqref{BB9}.
\end{proof}

We have the following useful corollary.
\bc(Khasminskii's type estimate)\label{Cor1}
In the same framework of Theorem \ref{Th41},  letting $\beta:=\frac{1}{2d+1}-\frac{1}{p}$ and $C$ be the same as in \eqref{BB9}, we have
\begin{enumerate}[(i)]
\item For each $m\in\mN$ and $0\leq t_0<t_1\leq T$, it holds that
\begin{align*}
\tilde\mE^{\tilde\sF_{t_0}}\left(\int^{t_1}_{t_0}f_s(\tilde Z_s)\dif s\right)^m\leq m!(C\|f\|_{\mL^p(T)}(t_1-t_0)^\beta)^m,
\end{align*}
where $\tilde\mE^{\tilde\sF_{t_0}}(\cdot)=\tilde\mE(\cdot|^{\tilde\sF_{t_0}})$.
\item For any $\lambda>0$ and $0\leq t_0<t_1\leq T$, it holds that
\begin{align*}
\tilde\mE^{\tilde\sF_{t_0}}\exp\left(\lambda\int^{t_1}_{t_0}f_s(\tilde Z_s)\dif s\right)\leq 2^{T (2C\lambda\|f\|_{\mL^p(T)})^{1/\beta}}.
\end{align*}
\end{enumerate}
\ec
\begin{proof}
(i) Still we drop the tilde below. For $m\in\mN$, noticing that
$$
\left(\int^{t_1}_{t_0}g(s)\dif s\right)^m=m!\int\!\!\!\cdots\!\!\!\int_{\Delta^m}g(s_1)\cdots g(s_m)\dif s_1\cdots\dif s_m,
$$
where
$$
\Delta^m:=\Big\{(s_1,\cdots, s_m):  t_0\leq s_1\leq s_2\leq\cdots\leq s_m\leq t_1\Big\},
$$
by \eqref{BB9}, we have
\begin{align*}
\mE^{\sF_{t_0}}\left(\int^{t_1}_{t_0}f_s(Z_s)\dif s\right)^m
&=m!\mE^{\sF_{t_0}}\left(\int\!\!\!\cdots\!\!\!\int_{\Delta^m}f_{s_1}(Z_{s_1})\cdots f_{s_m}(Z_{s_m})\dif s_1\cdots\dif s_m\right)\\
&=m!\mE^{\sF_{t_0}}\Bigg(\int\!\!\!\cdots\!\!\!\int_{\Delta^{m-1}}f_{s_1}(Z_{s_1})\cdots f_{s_{m-1}}(Z_{s_{m-1}})\\
&\quad\times\mE\left( \int^{t_1}_{s_{m-1}}f_{s_m}(Z_{s_m})\dif s_m\Big|_{\sF_{s_{m-1}}}\right)\dif s_1\cdots\dif s_{m-1}\Bigg)\\
&\leq m!\mE^{\sF_{t_0}}\int\!\!\!\cdots\!\!\!\int_{\Delta^{m-1}}f_{s_1}(Z_{s_1})\cdots f_{s_{m-1}}(Z_{s_{m-1}})\\
&\quad \times C(t_1-s_{m-1})^\beta\|f\|_{\mL^p(T)}\dif s_1\cdots\dif s_{m-1}\leq\cdots\\
&\leq m!(C\|f\|_{\mL^p(T)}(t_1-t_0)^\beta)^m.
\end{align*}
(ii) For $\lambda>0$, let us choose $n$ such that for $s_j=t_0+\frac{j(t_1-t_0)}{n}$,
$$
\lambda C\|f\|_{\mL^p(T)}(s_{j+1}-s_{j})^\beta\leq1/2.
$$
Then by (i) we have
\begin{align*}
\mE^{\sF_{s_j}}\exp\left(\lambda\int^{s_{j+1}}_{s_j}f_s(Z_s)\dif s\right)
&=\sum_{m}\frac{1}{m!}\mE^{\sF_{s_j}}\left(\lambda\int^{s_{j+1}}_{s_j}f_s(Z_s)\dif s\right)^m\leq 2.
\end{align*}
Hence,
\begin{align*}
&\mE^{\sF_{t_0}}\exp\left(\lambda\int^{t_1}_{t_0}f_s(Z_s)\dif s\right)=\mE^{\sF_{t_0}}\left(\prod_{j=0}^{n-1}\exp\left(\lambda\int^{s_{j+1}}_{s_j}f_s(Z_s)\dif s\right)\right)\\
&=\mE^{\sF_{t_0}}\left(\prod_{j=0}^{n-2}\exp\left(\lambda\int^{s_{j+1}}_{s_j}f_s(Z_s)\dif s\right)
\mE^{\sF_{s_{n-1}}}\exp\left(\lambda\int^{s_{n}}_{s_{n-1}}f_s(Z_s)\dif s\right)\right)\\
&\leq 2\mE^{\sF_{t_0}}\left(\prod_{j=0}^{n-2}\exp\left(\lambda\int^{s_{j+1}}_{s_j}f_s(Z_s)\dif s\right)\right)\leq\cdots\leq 2^{n}.
\end{align*}
The proof is complete.
\end{proof}
\subsection{Well-posedness of martingale problem} In this subsection we show the well-posedness of martingale problem for $\sL^{a,b}_t$. More precisely,
\bt\label{Th62}
Suppose that {\bf (UE)} holds, and for any $T>0$,
\begin{align}\label{CC9}
\lim_{|z-z'|\to 0}\sup_{t\in[0,T]}\|\sigma_t(z)-\sigma_t(z')\|=0,
\end{align}
and $b\in\mL^q(T)$ for some $q\in(2(2d+1),\infty]$. For each $(r,z)\in\mR_+\times\mR^{2d}$, 
the set $\sP^{\sigma,b}_{r,z}$ has one and only one point. 
In particular, the martingale problem for $\sL^{a,b}_t$ is well-posed.
\et
\begin{proof}
Below, we shall fix starting point $(r,z)\in\mR_+\times\mR^{2d}$ and divide the proof into three steps.

{\bf (a)} We first show the uniqueness. 
For $\varphi\in C^\infty_c(\mR^{2d})$ and $t_1>r$, let $u=\cR^{a,b}_{0,t_1}(\varphi)\in\mH^{2/3,2}_p(t_1)$, which satisfies
\begin{align}\label{KJ4}
\p_t u+\sL^{a,b}_t u+\varphi=0,\ \ u_{t_1}=0.
\end{align}
Let $u^\eps_t(z):=u_t*\varrho_\eps(z)$ and $\mP\in\sP^{\sigma,b}_{r,z}$.
By \eqref{MM4}, \eqref{BB9} and a standard approximation for the time variable, one sees that
$$
t\mapsto u^\eps(t,Z_{t})-\int^{t}_r(\p_s+\sL^{a,b}_s)u^\eps(s,Z_s)\dif s
$$ 
is an $\sF_t$-martingale under $\mP$ after time $r$. Thus, by \eqref{KJ4}, we have
\begin{align*}
u^\eps_r(z)&=-\mE\left(\int^{t_1}_r(\p_s u+\sL^{a,b}_s)u^\eps_s(Z_s)\dif s\right)\\
&=\mE\left(\int^{t_1}_r\big([\varrho_\eps, \sL^{a,b}_s]u_s+\varphi^\eps\big)(Z_s)\dif s\right).
\end{align*}
By Theorem \ref{Th41} and Lemma \ref{Le24}, taking limits $\eps\to 0$ for both sides yields
$$
u_r(z)=\mE\left(\int^{t_1}_r\varphi(Z_s)\dif s\right),\ t_1\geq r.
$$
In particular, we have  for any $\mP_1,\mP_2\in\sP^{\sigma,b}_{r,z}$ and $t\geq r$,
$$
\mE_{\mP_1}\varphi(Z_t)=\mE_{\mP_2}\varphi(Z_t).
$$
By \cite[Theorem 6.2.3]{St-Va}, we get the uniqueness.

{\bf (b)} For $n\in\mN$, let $\varrho_{1/n}$ be defined by \eqref{Rho} with $\eps=1/n$, and define
$$
b^n_t:=b_t*\varrho_{1/n}, \ \sigma^n_t:=\sigma_t*\varrho_{1/n}.
$$
Clearly,
$$
b^n\in L^q([0,T]; C^\infty_b(\mR^{2d})), \ \sigma^n\in L^\infty([0,T]; C^\infty_b(\mR^{2d}))
$$
and
$$
\|b^n\|_{\mL^q(T)}\leq\|b\|_{\mL^q(T)},\ \omega_{a^n}(\delta)\leq\omega_{a}(\delta).
$$  
By the classical theory of SDEs, the following SDE admits a unique strong solution $Z^n_t=(X^n_t,\dot X^n_t)$
$$
\dif Z^n_t=(\dot X^n_t,b^n_t(Z^n_t))\dif t+(0,\sigma^n_t(Z^n_t)\dif W_t),\ Z^n_t|_{[0,r]}=z.
$$
By Corollary \ref{Cor1}, for any $m\in\mN$ and $T>0$, there is a constant $C=C(m,T)>0$ such that for all $r\leq t_0<t_1\leq T$ and $n\in\mN$,
\begin{align*}
\mE\left|\int^{t_1}_{t_0}b^n_s(Z^n_s)\dif s\right|^m\leq C\|b\|_{\mL^q(T)}^m(t_1-t_0)^{m(\frac{1}{2d+1}-\frac{1}{p})}.
\end{align*}
Let $\mP_n$  be the probability distribution of $Z^n$ in $\Omega$. By the above moment estimate, 
it is by now standard to show that $(\mP_n)_{n\in\mN}$ is tight.  
 
 {\bf (c)} By extracting a subsequence, without loss of generality, we may assume that
 $\mP_n$ weakly converges to some probability measure $\mP$.  To see that $\mP\in\sP^{\sigma,b}_{r,z}$,
it suffices to show that for any $\varphi\in C^\infty_c(\mR^{2d})$, $M^{r,\varphi}_t$ defined by \eqref{MM4} is an $\sF_t$-martingale
under $\mP$. Equivalently, for any $r\leq t_0\leq t_1$ and any bounded continuous $\sF_{t_0}$-measurable $G$,
\begin{align}\label{MM5}
\mE_{\mP}(M^{r,\varphi}_{t_1}G)=\mE_{\mP}(M^{r,\varphi}_{t_0}G).
\end{align}
Notice that
\begin{align}\label{Lim3}
\mE_{\mP_n}(M^{r,\varphi}_{n,t_1}G)=\mE_{\mP_n}(M^{r,\varphi}_{n,t_0}G),
\end{align}
where 
$$
M^{r,\varphi}_{n,t_i}:=\varphi(Z_{t_i})-\int^{t_i}_r\sL^{a^n,b^n}_s \varphi(Z_s)\dif s,\ i=0,1.
$$
Let us prove the following limit: for $i=0,1$,
\begin{align}\label{Lim0}
\lim_{n\to\infty}\mE_{\mP_n}\left(G\int^{t_i}_r(b^n_s\cdot\nabla_\v\varphi)(Z_s)\dif s\right)=\mE_{\mP}\left(G\int^{t_i}_r(b_s\cdot\nabla_\v\varphi)(Z_s)\dif s\right).
\end{align}
For any $p\in(2d+1,q)$ and $T>0$, by Theorem \ref{Th41}, there is a constant $C>0$ such that for all $n\in\mN$ and $f\in \mL^p(T)$,
\begin{align*}
\mE_{\mP_n}\left(\int^T_r f_s(Z_s)\dif s\right)\leq C\|f\|_{\mL^p(T)}.
\end{align*} 
Let $f\in C_c([0,T]\times\mR^{2d})$. By taking weak limits, we have 
\begin{align*}
\mE_{\mP}\left(\int^T_r f_s(Z_s)\dif s\right)&=\lim_{n\to\infty}\mE_{\mP_n}\left(\int^T_r f_s(Z_s)\dif s\right)\leq C\|f\|_{\mL^p(T)}.
\end{align*}
By a monotone class argument, the above estimate still holds for all $f\in\mL^p(T)$. 
Let the support of $\varphi$ be contained in $B_R$. Thus, if we let $\mP_\infty=\mP$, then
\begin{align}\label{Lim1}
\begin{split}
&\lim_{m\to\infty}\sup_{n\in\mN\cup\{\infty\}}\mE_{\mP_n}\left(G\int^{t_i}_r|(b^m_s-b_s)\cdot\nabla_\v\varphi|(Z_s)\dif s\right)\\
&\quad\leq C\|G\|_\infty\|\nabla_\v\varphi\|_\infty\lim_{m\to\infty}\|(b^m-b)1_{B_R}\|_{\mL^p(r,t_i)}=0.
\end{split}
\end{align}
On the other hand, for each $m\in\mN$, since $\omega\mapsto G(\omega)\int^{t_i}_r(b^m_s\cdot\nabla_\v\varphi)(Z_s(\omega))\dif s$ 
is a continuous and bounded functional, we have
$$
\lim_{n\to\infty}\mE_{\mP_n}\left(G\int^{t_i}_r(b^m_s\cdot\nabla_\v\varphi)(Z_s)\dif s\right)=\mE_{\mP}\left(G\int^{t_i}_r(b^m_s\cdot\nabla_\v\varphi)(Z_s)\dif s\right).
$$
Combining this with \eqref{Lim1}, we get \eqref{Lim0}. Similarly, one can show
\begin{align}\label{Lim2}
\lim_{n\to\infty}\mE_{\mP_n}\left(G\int^{t_i}_r\sL^{a^n,0}_r\varphi(Z_s)\dif s\right)=\mE_{\mP}\left(G\int^{t_i}_r\sL^{a,0}_r\varphi(Z_s)\dif s\right).
\end{align}
Finally, by taking weak limits for both sides of \eqref{Lim3} and using \eqref{Lim0} and \eqref{Lim2}, we get \eqref{MM5}. The proof is complete.
\end{proof}

\br
When $b$ is bounded measurable ($q=\infty$), this result has been proven in \cite{Pr} and \cite{Me}. Therein, more general equations were considered.
However, by comparing with the original proof of Stroock and Varadhan \cite[Chapter 7]{St-Va},
our proof is quite different from \cite{Pr, Me} as our starting point is
a global apriori Krylov's estimate (see Theorem \ref{Th41}). In principle, it is reasonable to believe 
that our argument is applicable for more general equations as studied in \cite{Pr, Me}.
\er
\br
By suitable localization techniques as developed in \cite{Pr}, it is possible to weaken the global assumptions in Theorem \ref{Th62} as local ones together with
some non-explosion conditions. 
\er

\subsection{Proof of Theorem \ref{Main2}}
Let $\nu\in\cP(\mR^{2d})$. By Theorem \ref{Th62}, the probability measure $\mP_\nu(A):=\int_{\mR^{2d}}\mP_{0,z}(A)\nu(\dif z)$ is the unique martingale
solution for $\sL^{a,b}_t$ starting from $\nu$ at time $0$. The conclusion of Theorem \ref{Main2} now follows by \cite[Theorem 2.5]{Tre}.

\section{Proof of Theorem \ref{Main}}

In this section we assume that $\sigma$ satisfies {\bf (UE)} and for some $p>2(2d+1)$,
\begin{align}\label{Con5}
\kappa_0:=\|b\|_{L^p(\mR_+; H^{2/3,0}_p)}+\|\nabla\sigma\|_{L^\infty(\mR_+; L^p)}<\infty.
\end{align}
For $n\in\mN$, let $\varrho_{1/n}$ be defined by \eqref{Rho} with $\eps=1/n$, as in the previous,  define
\begin{align}\label{Reg}
b^n_t:=b_t*\varrho_{1/n}, \ \sigma^n_t:=\sigma_t*\varrho_{1/n},\ \ b^\infty_t:=b_t, \ \sigma^\infty_t:=\sigma_t.
\end{align}
\bl\label{Le51}
Assume {\bf (UE)} and \eqref{Con5}. Then {\bf (H$^{\delta,p}_K$)} holds for $a^n:=\frac{1}{2}\sigma^n(\sigma^n)^*$ 
uniformly with respect to $n$, and there exists a constant $C=C(d,p)>0$ such that for all $s$,
$$
\|\sigma^n_s-\sigma_s\|_\infty\leq C\|\nabla \sigma_s\|_p n^{\frac{2d}{p}-1},\ \ \sup_{\v}\|\Delta_x^{\frac{1}{3}}\sigma^n_s(\cdot,\v)\|_p\leq C\|\nabla\sigma_s\|_p.
$$
\el
\begin{proof}
Since $p>2d$, by Morrey's inequality, there is a constant $C=C(d,p)>0$ such that
$$
|\sigma_s(z)-\sigma_s(z')|\leq C|z-z'|^{1-\frac{2d}{p}}\|\nabla \sigma_s\|_p,\ \ z,z'\in\mR^{2d}.
$$
From this, it is easy to see that {\bf (H$^{\delta,p}_K$)} holds for $a^n=\frac{1}{2}\sigma^n(\sigma^n)^*$  uniformly with respect to $n$, and
\begin{align*}
\|\sigma^n_s-\sigma_s\|_\infty&\leq\int_{\mR^{2d}}\|\sigma_s(\cdot+z)-\sigma_s(\cdot)\|_\infty\varrho_{1/n}(z)\dif z\\
&\leq C\|\nabla \sigma_s\|_p\int_{\mR^{2d}}|z|^{1-\frac{2d}{p}}\varrho_{1/n}(z)\dif z\leq C\|\nabla \sigma_s\|_p n^{\frac{2d}{p}-1}.
\end{align*}
Moreover, since $p>4d$, by \eqref{Em0} and \eqref{Em88}, we have
\begin{align*}
\sup_{\v}\|\Delta_x^{\frac{1}{3}}\sigma^n_s(\cdot,\v)\|^p_p
&\leq\int_{\mR^d}\sup_\v |\Delta_x^{\frac{1}{3}}\sigma^n_s(x,\v)|^p_p\dif x
\preceq\|\Delta_x^{\frac{1}{3}}\sigma^n_s\|^p_p+\|\Delta_x^{\frac{1}{3}}\Delta_\v^{\frac{1}{8}}\sigma^n_s\|^p_{p}\\
&\preceq\|\Delta^{\frac{11}{24}}\sigma^n_s\|^p_{p}
\preceq\|\Delta^{\frac{1}{2}}\sigma^n_s\|^p_{p}\stackrel{\eqref{Riesz}}{\preceq}\|\nabla\sigma^n_s\|^p_p.
\end{align*}
Here $\Delta=\Delta_x+\Delta_\v$. The proof is complete.
\end{proof}

Let $T>0$ and $\lambda\geq 1$. For $n\in\mN\cup\{\infty\}$, let $\u^n_\lambda\in \mH^{0,2}_p(T)$ uniquely solve the following PDE:
$$
\p_t\u^n_\lambda+\sL^{a^n,b^n}_t\u^n_\lambda-\lambda \u^n_\lambda+b^n=0,\ \ \u^n_{_\lambda,T}=0,
$$
where $a^n:=\frac{1}{2}\sigma^n(\sigma^n)^*$.
By Lemma \ref{Le51} and Theorem \ref{Th32}, there is a constant $C=C(d,p,\kappa_0,K)>0$
such that for all $n\in\mN\cup\{\infty\}$,
\begin{align}\label{PP1}
\|\nabla\nabla_\v\u^n_\lambda\|_{\mL^p(T)}\leq C\|b\|_{\mH^{2/3,0}_p(T)},
\end{align}
and by \eqref{E28}, \eqref{E29},  \eqref{Pri00} with $\alpha=\frac{4d+1}{3(2d+1)}$ and $\beta=\frac{4d+1}{2d+1}$,
\begin{align}
\|\nabla_\v\u^n_\lambda\|_{\mL^\infty(T)}\leq C\|\u^n_\lambda\|_{L^\infty(0,T; \mH^{\alpha,\beta}_p)}
\leq C(T\wedge\lambda^{-1})^{\frac{1}{2(2d+1)}-\frac{1}{p}}\|b\|_{\mL^p(T)}.\label{PP2}
\end{align}
\bl\label{Le52}
Under {\bf (UE)} and \eqref{Con5}, there is a constant  $C>0$ depending only on $d,p,\kappa_0,K, T,\lambda$ such that for all $n,m\in\mN$,
$$
\|\u^n_\lambda-\u^m_\lambda\|_{\mL^\infty(T)}+\|\nabla_\v\u^n_\lambda-\nabla_\v\u^m_\lambda\|_{\mL^\infty(T)}
\leq C\Big(\|b^n-b^m\|_{\mL^p(T)}+(n\wedge m)^{\frac{2d}{p}-1}\Big).
$$
\el
\begin{proof}
Let $\w^{n,m}:=\u^n_\lambda-\u^m_\lambda$. By equation \eqref{PDE}, we have
$$
\p_t\w^{n,m}+(\sL^{a^n,b^n}_t-\lambda)\w^{n,m}+(\sL^{a^n,b^n}_t-\sL^{a^m,b^m}_t)\u^m_\lambda+b^n-b^m=0.
$$
Noticing that
$$
g^{n,m}_{\lambda,t}:=(\sL^{a^n,b^n}_t-\sL^{\sigma^m,b^m}_t)\u^m_{\lambda,t}=\tr((a^n_t-a^m_t)\cdot\nabla^2_\v \u^m_{\lambda,t})
+(b^n_t-b^m_t)\cdot\nabla_\v \u^m_{\lambda,t},
$$
by \eqref{PP1}, \eqref{PP2} and Lemma \ref{Le51}, we have
\begin{align*}
\|g^{n,m}_\lambda\|_{\mL^p(T)}&\preceq\|a^n-a^m\|_{\mL^\infty(T)}\|\nabla^2_\v \u^m_\lambda\|_{\mL^p(T)}
+\|b^n-b^m\|_{\mL^p(T)}\|\nabla_\v \u^m_\lambda\|_{\mL^\infty(T)}\\
&\preceq (n\wedge m)^{\frac{2d}{p}-1}+\|b^n-b^m\|_{\mL^p(T)}.
\end{align*}
By \eqref{E29} and \eqref{Pri00} with $\alpha=\frac{4d}{3(2d+1)}$ and $\beta=\frac{4d}{2d+1}$, we have
\begin{align*}
\|\w^{n,m}_t\|_\infty\preceq\|\w^{n,m}_t\|_{\alpha,\beta;p}&\preceq \|g^{n,m}_\lambda\|_{\mL^p(T)}+\|b^n-b^m\|_{\mL^p(T)}\\
&\preceq (n\wedge m)^{\frac{2d}{p}-1}+\|b^n-b^m\|_{\mL^p(T)},
\end{align*}
and by \eqref{E28}, \eqref{E29},  \eqref{Pri00} with $\alpha=\frac{4d+1}{3(2d+1)}$ and $\beta=\frac{4d+1}{2d+1}$,
\begin{align*}
\|\nabla_\v\w^{n,m}_t\|_\infty&\preceq \|\nabla_\v\w^{n,m}_t\|_{\alpha(\beta-1)/\beta,\beta-1;p}
\preceq \|\w^{n,m}_t\|_{\alpha,\beta;p}\\
&\preceq (n\wedge m)^{\frac{2d}{p}-1}+\|b^n-b^m\|_{\mL^p(T)}.
\end{align*}
The proof is complete.
\end{proof}
For $n\in\mN\cup\{\infty\}$, let 
\begin{align}\label{HH}
H^n_t(x,\v):=\v+\u^n_{\lambda,t}(x,\v).
\end{align}
By \eqref{PP2}, one can choose $\lambda$ large enough (being independent of $n$ and fixed below) so that
\begin{align}\label{T0}
\|\nabla_\v\u^n_\lambda\|_{\mL^\infty(T)}\leq\tfrac{1}{2},
\end{align}
and thus,
\begin{align}\label{DIF}
\tfrac{1}{2}|\v-\v'|\leq |H^n_t(x,\v)-H^n_t(x,\v')|\leq \tfrac{3}{2}|\v-\v'|.
\end{align}
Observing that
$$
\p_tH^n+\sL^{a^n,b^n}_tH^n-\lambda \u^n_\lambda=0,
$$
by It\^o's formula, we have
\begin{align}\label{SDE1}
H^n_t(Z^n_t)=H^n_0(Z^n_0)+\lambda\int^t_0\u^n_{\lambda,s}(Z^n_s)\dif s+\int^t_0\Theta^n_s(Z^n_s)\dif W_s,
\end{align}
where $\Theta^n_s(z):=(\nabla_\v H^n_s\cdot \sigma^n_s)(z)$ satisfies by \eqref{T0} and \eqref{PP1} that
$$
\|\Theta^n\|_\infty\leq 2\|\sigma\|_\infty,
$$
and for the given $p>2(2d+1)$,
\begin{align}\label{633}
\begin{split}
\|\nabla\Theta^n\|_{\mL^p(T_0)}&\leq 2\|\nabla\sigma^n\|_{\mL^p(T_0)}+\|\sigma^n\|_\infty\|\nabla\nabla_\v H^n\|_{\mL^p(T_0)}\\
&\leq 2\|\nabla\sigma\|_{\mL^p(T_0)}+C\|\sigma\|_\infty\|b\|_{\mH^{2/3,0}_p(T_0)}.
\end{split}
\end{align}
For $n\in\mN$, since $b^n\in L^p([0,T]; C^\infty_b(\mR^{2d}))$ and
$\sigma^n\in L^\infty([0,T]; C^\infty_b(\mR^{2d}))$,  the following SDE admits a unique solution $Z^n_t=(X^n_t,\dot X^n_t)$
\begin{align}\label{App0}
\dif Z^n_t=(\dot X^n_t,b^n_t(Z^n_t))\dif t+(0,\sigma^n_t(Z^n_t)\dif W_t),\ \ Z^n_0=z=(x,\v)\in\mR^{2d}.
\end{align}
We have
\bl\label{Le53}
Under {\bf (UE)} and \eqref{Con5}, for any $q\geq 2$, there is a constant $C>0$ such that for all $n,m\in\mN$,
\begin{align}\label{PP4}
\left\|\sup_{t\in[0,T_0]}|Z^n_t-Z^m_t|\right\|_{L^q(\Omega)}\leq 
C\left(\|b^n-b^m\|_{\mL^p(T_0)}+(n\wedge m)^{\frac{2d}{p}-1}\right).
\end{align}
\el
\begin{proof}
By \eqref{SDE1} and It\^o's formula, we have
\begin{align*}
&|H^n_t(Z^n_t)-H^m_t(Z^m_t)|^2=|H^n_0(z)-H^m_0(z)|^2+\int^t_0\|\Theta^n_s(Z^n_s)-\Theta^m_s(Z^m_s)\|^2\dif s\\
&\qquad\qquad+2\lambda\int^t_0\<H^n_s(Z^n_s)-H^m_s(Z^m_s), \u^n_{\lambda,s}(Z^n_s)-\u^m_{\lambda,s}(Z^m_s)\>\dif s\\
&\qquad\qquad+2\int^t_0\<H^n_s(Z^n_s)-H^m_s(Z^m_s), (\Theta^n_s(Z^n_s)-\Theta^m_s(Z^m_s))\dif W_s\>,
\end{align*}
and also,
$$
|X^n_t-X^m_t|^2=2\int^t_0\<X^n_s-X^m_s,\dot X^n_s-\dot X^m_s\>\dif s.
$$
If we set
$$
\xi_t:=|H^n_t(Z^n_t)-H^n_t(Z^m_t)|^2+|X^n_t-X^m_t|^2,
$$
then 
\begin{align*}
\xi_t&\leq 2\|H^n_t-H^m_t\|^2_\infty+2|H^n_t(Z^n_t)-H^m_t(Z^m_t)|^2+|X^n_t-X^m_t|^2\\
&\leq\xi_0+\int^t_0\zeta^{(1)}_s\dif s+\int^t_0\zeta^{(2)}_s\dif W_s+\int^t_0\xi_s\beta_s\dif s+\int^t_0\xi_s\alpha_s\dif W_s,
\end{align*}
where $\xi_0:=3\|\u^n-\u^m\|^2_{\mL^\infty(T_0)}$ and
\begin{align*}
\zeta^{(1)}_s&:=4\|\Theta^n_s(Z^m_s)-\Theta^m_s(Z^m_s)\|^2+6\lambda|H^n_s(Z^m_s)-H^m_s(Z^m_s)|^2,\\
\zeta^{(2)}_s&:=4(\Theta^n_s(Z^n_s)-\Theta^m_s(Z^m_s))^*(H^n_s(Z^n_s)-H^m_s(Z^m_s))\\
&\quad-4(\Theta^n_s(Z^n_s)-\Theta^n_s(Z^m_s))^*(H^n_s(Z^n_s)-H^n_s(Z^m_s)),\\
\beta_s&:=4\|\Theta^n_s(Z^n_s)-\Theta^n_s(Z^m_s)\|^2/\xi_s+6\lambda+\lambda|\dot X^n_s-\dot X^m_s|^2/\xi_s\\
&\quad+2\<X^n_s-X^m_s,\dot X^n_s-\dot X^m_s\>/\xi_s,\\
\alpha_s&:=4(\Theta^n_s(Z^n_s)-\Theta^n_s(Z^m_s))^*(H^n_s(Z^n_s)-H^n_s(Z^m_s))/\xi_s.
\end{align*}
By \eqref{DIF}, one has
\begin{align}\label{Es31}
\xi_t\asymp|Z^n_t-Z^m_t|^2=|X^n_t-X^m_t|^2+|\dot X^n_t-\dot X^m_t|^2,
\end{align}
which implies by \eqref{Es2} that
\begin{align}\label{PP3}
|\beta_s|+|\alpha_s|^2\preceq 1+\big(\cM|\nabla \Theta^n_s|(Z^n_s)\big)^2+\big(\cM|\nabla \Theta^n_s|(Z^m_s)\big)^2=:G^{n,m}_s.
\end{align}
In view of $p>2(2d+1)$ and $\sup_n\|\cM|\nabla \Theta^n|\|_{\mL^p(T)}<\infty$, by (ii) of Corollary \ref{Cor1}, we have for any $\gamma>0$,
$$
\sup_{n,m}\mE\exp\left\{\gamma\int^{T}_0G^{n,m}_s\dif s\right\}<\infty.
$$
Therefore, by Lemma \ref{Le35} below with $q_0=q$, $q_1=q_2=q_3=3q/2$, we have
\begin{align}\label{513}
\|\xi^*_{T}\|_q\preceq \|\u^n-\u^m\|^2_{\mL^\infty(T)}+\left\|\int^{T}_0|\zeta^{(1)}_s|\dif s\right\|_{3q/2}
+\left\|\int^{T}_0|\zeta^{(2)}_s|^2\dif s\right\|^{1/2}_{3q/4},
\end{align}
where $\xi^*_{T}:=\sup_{t\in[0,T]}\xi_t$.

Now, letting
$$
\ell_{n,m}:=\|\sigma^n-\sigma^m\|_{\mL^\infty(T)}+\|\u^n-\u^m\|_{\mL^\infty(T)}+\|\nabla\u^n-\nabla\u^m\|_{\mL^\infty(T)},
$$
by definitions, we have
\begin{align*}
|\zeta_s^{(1)}|\preceq\ell^2_{n,m},
\end{align*}
and by \eqref{Es2} and \eqref{Es31},
\begin{align*}
|\zeta^{(2)}_s|^2&\preceq |(\Theta^n_s(Z^n_s)-\Theta^m_s(Z^m_s))^*(H^n_s(Z^m_s)-H^m_s(Z^m_s))|^2\\
&\quad+|(\Theta^n_s(Z^m_s)-\Theta^m_s(Z^m_s))^*(H^n_s(Z^n_s)-H^n_s(Z^m_s))|^2\\
&\preceq \Big(\|\Theta^n_s(Z^m_s)-\Theta^m_s(Z^m_s)\|^2+\|\Theta^n_s(Z^n_s)-\Theta^n_s(Z^m_s)\|^2\Big)\\
&\quad\times|\u^n_s(Z^m_s)-\u^m_s(Z^m_s)|^2+\|\Theta^n_s(Z^m_s)-\Theta^m_s(Z^m_s)\|^2\xi_s\\
&\preceq \ell^4_{n,m}+\ell^2_{n,m}G^{n,m}_s\xi_s,
\end{align*}
where $G^{n,m}_s$ is defined by \eqref{PP3}.
Substituting these estimates into \eqref{513} and by H\"older's inequlity and Young's inequality, 
we have for any $\eps\in(0,1)$,
\begin{align}
\|\xi^*_{T}\|_q
&\preceq\ell^2_{n,m}+\ell_{n,m}\left\|\int^{T}_0G^{n,m}_s\xi_s\dif s\right\|^{1/2}_{3q/4}\no\\
&\preceq\ell^2_{n,m}+\ell_{n,m}\left\|\xi^*_{T}\int^{T}_0G^{n,m}_s\dif s\right\|^{1/2}_{3q/4}\no\\
&\preceq\ell^2_{n,m}+\ell_{n,m}\left\|\int^{T}_0G^{n,m}_s\dif s\right\|^{1/2}_{3q}\|\xi^*_{T}\|_q^{\frac{1}{2}}\no\\
&\preceq C_\eps\ell^2_{n,m}+\eps\|\xi^*_{T}\|_q.\label{ERQ}
\end{align}
Letting $\eps$ be small enough and noticing that by Lemmas \ref{Le51} and \ref{Le52},
$$
\ell_{n,m}\leq C\left(\|b^n-b^m\|_{\mL^p(T)}+(n\wedge m)^{\frac{2d}{p}-1}\right),
$$
we obtain the desired estimate \eqref{PP4} from \eqref{ERQ}.
\end{proof}
\bt(Existence of strong solutions)\label{Th54}
Under {\bf (UE)} and \eqref{Con5}, for any $q\geq 2$, 
there is a continuous $\sF_t$-adapted process $Z_t=(X_t,\dot X_t)$ solving equation \eqref{SDE0}, and such that
\begin{align}\label{PP7}
\left\|\sup_{t\in[0,T]}|Z^n_t-Z_t|\right\|_{L^q(\Omega)}\leq C\left(\|b^n-b\|_{\mL^p(T)}+n^{\frac{2d}{p}-1}\right).
\end{align}
\et
\begin{proof}
First of all, by \eqref{PP4}, there exists a continuous $\sF_t$-adapted process $Z^\infty_t:=Z_t=(X_t,\dot X_t)$ satisfying \eqref{PP7}. 
Moreover, for any $p>2d+1$, there is a positive constant $C=C(T,K,\kappa_0,d,p)$ such that for any $f\in\mL^p(T)$, 
\begin{align}\label{Est1}
\mE\left(\int^{T}_0f_s(Z^\infty_s)\dif s\right)\leq C\|f\|_{\mL^p(T)}.
\end{align}
In fact,  for any $f\in C_c([0,T]\times\mR^{2d})$, by \eqref{BB9} we have
\begin{align}\label{Est2}
\mE\left(\int^{T}_0f_s(Z^n_s)\dif s\right)\leq C\|f\|_{\mL^p(T)},
\end{align}
which gives \eqref{Est1} by taking limit $n\to\infty$ and a monotone class argument.

Below we show that $Z_t$ solves SDE \eqref{SDE0} by taking limits for \eqref{App0}. For this, we only need to show  the following two limits:
\begin{align}
&\lim_{n\to\infty}\mE\left|\int^t_0 b^n_s(Z^n_s)\dif s-\int^t_0 b_s(Z_s)\dif s\right|=0,\label{Lim7}\\
&\lim_{n\to\infty}\mE\left|\int^t_0 \sigma^n_s(Z^n_s)\dif W_s-\int^t_0 \sigma_s(Z_s)\dif W_s\right|=0.\label{Lim8}
\end{align}
For \eqref{Lim7}, by \eqref{Est1} and \eqref{Est2} we have
$$
\lim_{n\to\infty}\sup_{m\in\mN\cup\{\infty\}}\mE\left(\int^t_0 |b^n_s(Z^m_s)-b_s(Z^m_s)|\dif s\right)
\leq C\lim_{n\to\infty}\|b^n-b\|_{\mL^p(t)}=0,
$$
and for each $n\in\mN$, by the dominated convergence theorem,
$$
\lim_{m\to\infty}\mE\left(\int^t_0 |b^n_s(Z^m_s)-b^n_s(Z^\infty_s)|\dif s\right)=0.
$$
Limit \eqref{Lim8} is similar. The proof is complete.
\end{proof}
Now, by a result of Cherny \cite{Che},  the existence of strong solutions together with the weak uniqueness (see Theorem \ref{Th62}) implies the pathwise uniqueness.
However, to show the homeomorphism property of $z\mapsto Z_t(z)$, 
one needs the following $q$-order moment estimate for all $q\in\mR$, which clearly implies the pathwise uniqueness.
\bl\label{Le55}
For any $q\in\mR$, there is a constant $C>0$ such that for any two solutions $Z_t$ and $Z_t'$ of SDE \eqref{SDE0}
with starting points $z=(x,\v)$ and $z'=(x',\v')$ respectively,
\begin{align}\label{Es7}
\mE\left(\sup_{t\in[0,T]}|Z_t-Z_t'|^{2q}\right)\leq C|z-z'|^{2q}.
\end{align}
Moreover, we also have
\begin{align}\label{Es8}
\mE\left(\sup_{t\in[0,T]}(1+|Z_t|^2)^{q}\right)\leq C(1+|z|^2)^{q}.
\end{align}
\el
\begin{proof}
Let $H=H^\infty$ be defined by \eqref{HH}. By \eqref{App0} and It\^o's formula, we have
\begin{align*}
|H_t(Z_t)-H_t(Z'_t)|^2&=|H_0(Z_0)-H_0(Z'_0)|^2+\int^t_0\|\Theta_s(Z_s)-\Theta_s(Z'_s)\|^2\dif s\\
&+2\lambda\int^t_0\<H_s(Z_s)-H_s(Z'_s), \u_{\lambda,s}(Z_s)-\u_{\lambda,s}(Z'_s)\>\dif s\\
&+2\int^t_0\<H_s(Z_s)-H_s(Z'_s),(\Theta_s(Z_s)-\Theta_s(Z'_s))\dif W_s\>.
\end{align*}
As in the proof of Lemma \ref{Le53}, if we set
$$
\xi_t:=|H_t(Z_t)-H_t(Z'_t)|^2+|X_t-X'_t|^2,
$$
then
$$
\xi_t=\xi_0+\int^t_0\xi_s\beta_s\dif s+\int^t_0\xi_s\alpha_s\dif W_s,
$$
where
\begin{align*}
\beta_s&:=[\|\Theta_s(Z_s)-\Theta_s(Z'_s)\|^2+2\<X_s-X'_s,\dot X_s-\dot X'_s\>]/\xi_s\\
&\quad+2\lambda\<H_s(Z_s)-H_s(Z'_s), \u_{\lambda,s}(Z_s)-\u_{\lambda,s}(Z'_s)\>/\xi_s
\end{align*}
and
$$
\alpha_s:=2(\Theta_s(Z_s)-\Theta_s(Z'_s))^*(H_s(Z_s)-H_s(Z'_s))/\xi_s.
$$
By Dol\`eans-Dade's exponential formula, we have
\begin{align}\label{H6}
\xi_t=\xi_0\exp\left\{\int^t_0\alpha_s\dif W_s+\int^t_0\Big[\beta_s-\tfrac{1}{2}|\alpha_s|^2\Big]\dif s\right\}
\end{align}
By \eqref{DIF}, one has
$$
\xi_t\asymp|Z_t-Z'_t|^2=|X_t-X'_t|^2+|\dot X_t-\dot X'_t|^2.
$$
Hence, by \eqref{HH} and \eqref{Es2} below,
$$
|\beta(s)|+|\alpha(s)|^2\preceq 1+(\cM|\nabla \Theta_s|(Z_s))^2+(\cM|\nabla \Theta_s|(Z'_s))^2.
$$
Since $\|\cM|\nabla \Theta|\|_{\mL^p(T)}<\infty$ and $p>2(2d+1)$, by Corollary \ref{Cor1}, we have for any $\gamma>0$,
\begin{align}\label{Exp}
\mE\exp\left\{\gamma\int^{T}_0\big(|\beta_s|+|\alpha_s|^2\big)\dif s\right\}<\infty.
\end{align}
For $q\in\mR$, let
$$
\sE^q_t:=\exp\left\{q\int^t_0\alpha_s\dif W_s-\tfrac{q^2}{2}\int^t_0|\alpha_s|^2\dif s\right\}.
$$
By \eqref{Exp} and Novikov's criterion, $t\mapsto \sE^q_t$ is a continuous martingale. Therefore, by \eqref{H6}, H\"older's inequality
and Doob's maximal inequality,
\begin{align*}
&\mE\left(\sup_{t\in[0,T]}|Z_t-Z'_t|^{2q}\right)\\
&\preceq|z-z'|^{2q}
\mE\left(\sup_{t\in[0,T]}\exp\left\{q\int^t_0\alpha_s\dif W_s+q\int^t_0[\beta_s-\tfrac{1}{2}|\alpha_s|^2]\dif s\right\}\right)\\
&\preceq|z-z'|^{2q}
\left(\mE\left((\sE^q)^*_{T}\right)^2\right)^{\frac{1}{2}}
\left(\mE\exp\left\{\int^{T}_0\big((q^2+q)|\alpha_s|^2+q\beta_s\big)\dif s\right\}\right)^{\frac{1}{2}}\\
&\preceq|z-z'|^{2q}\left(\mE(\sE^{q}_{T})^2\right)^{\frac{1}{2}}
\preceq|z-z'|^{2q},
\end{align*}
which gives \eqref{Es7}.

On the other hand, if we let 
$$
\xi_t:=1+|H_t(Z_t)|^2+|X_t|^2,
$$
then for any $q\in\mR$, by It\^o's formula, we have
\begin{align*}
&\xi_t^q=\xi_0^q+q\int^t_0\xi_s^{q-1}(2\<X_s,\dot X_s\>+\|\Theta_s(Z_s)\|^2-\lambda\<H_s(Z_s),\u_{\lambda,s}(Z_s)\>)\dif s\\
&+2q\int^t_0\xi_s^{q-1}\<H_s(Z_s), \Theta_s(Z_s)\dif W_s\>+2q(q-1)\int^t_0\xi_s^{q-2}|\Theta^*_s(Z_s)H_s(Z_s)|^2\dif s.
\end{align*}
Noticing that
$$
|H_t(Z_t)|\asymp1+|\dot X_t|,\ \ \xi_t\asymp 1+|Z_t|^2=1+|X_t|^2+|\dot X_t|^2,
$$
by Burkholder's inequality, we have
\begin{align*}
\mE\left(\sup_{s\in[0,t]}\xi_t^{2q}\right)\preceq \xi_0^{2q}+\mE\int^t_0\xi_s^{2q}\dif s,
\end{align*}
which in turn gives \eqref{Es8} by Gronwall's inequality.
\end{proof}

Now we can give
\begin{proof}[Proof of Theorem \ref{Main}]
The existence  and uniqueness of a strong solution in time interval $[0,T]$ 
follows by Theorem \ref{Th54} and Lemma \ref{Le55}. The bi-continuous version of $(t,z)\mapsto Z_t(z)$ 
follows by \eqref{Es7} and Kolmogorov's continuity criterion. 

{\bf (A)} As for the homeomorphism property, it follows by Lemma \ref{Le55} and 
Kunita's argument (see \cite{Ku, Zh2}). 

{\bf (B)} The weak differentiability of $z\mapsto Z_t(z)$ and estimate \eqref{YN1} follow by \eqref{Es7} 
and \cite[Theorem 1.1]{Xi-Zh}. 

{\bf (C)} It follows by \eqref{PP7}.
\end{proof}

\section{Appendix}

The following stochastic Gronwall's type lemma is probably well-known. Since we can not find it in the literature, a proof is provided here for the reader's convenience.
\bl\label{Le35}
For given $T>0$, let $(\xi_t)_{t\in [0, T]}$ and $(\beta_t)_{t\in [0, T]}$ (resp. $(\alpha_t)_{t\in [0, T]}$) be 
two real-valued (resp. $\mR^d$-valued) measurable $\sF_t$-adapted processes. Let $\zeta_t$ be an It\^o process with the form:
$$
\zeta_t=\zeta_0+\int^t_0\zeta^{(1)}_s\dif s+\int^t_0\zeta^{(2)}_s\dif W_s.
$$
Suppose that for any $\gamma>0$,
\begin{align}
\kappa_\gamma:=\mE\exp\left\{\gamma\int^T_0\big(|\beta_s|+|\alpha_s|^2\big)\dif s\right\}<\infty,\label{C1}
\end{align}
and 
\begin{align}\label{C2}
0\leq\xi_t\leq\zeta_t+\int^t_0\xi_s\beta_s\dif s+\int^t_0\xi_s\alpha_s\dif W_s.
\end{align}
Then for any $q_0\in[1,\infty)$ and $q_1,q_2,q_3>q_0$, there is a constant $C>0$ only depending on $q_i,\kappa_{\gamma}, i=0,1,2,3$ such that
\begin{align}
\|\xi^*_T\|_{q_0}\leq C\left(\|\zeta_0\|_{q_1}+\left\|\int^T_0|\zeta_s^{(1)}|\dif s\right\|_{q_2}
+\left\|\int^T_0|\zeta_s^{(2)}|^2\dif s\right\|^{1/2}_{q_3/2}\right),\label{HG8}
\end{align}
where $\xi^*_T:=\sup_{t\in[0,T]}\xi_t$ and $\|\cdot\|_{q_i}$ denotes the norm in $L^{q_i}(\Omega)$.
\el
\begin{proof}
Write
\begin{align*}
\eta_t&:=\zeta_t+\int^t_0\xi_s\beta_s\dif s+\int^t_0\xi_s\alpha_s\dif W_s\\
&=\zeta_t+\int^t_0\eta_s\bar\beta_s\dif s+\int^t_0\eta_s\bar\alpha_s\dif W_s,
\end{align*}
where $\bar\beta_s:=\xi_s\beta_s/\eta_s$ and $\bar\alpha_s:=\xi_s\alpha_s/\eta_s$. Here we use the convention $\frac{0}{0}:=0$.

Define
$$
M_t:=\exp\left\{\int^t_0\bar\alpha_s\dif W_s+\int^t_0(\bar\beta_s-\tfrac{1}{2}|\bar\alpha_s|^2)\dif s\right\}.
$$
By It\^o's formula, we have
$$
M_t=1+\int^t_0M_s\bar\beta_s\dif s+\int^t_0M_s\bar\alpha_s\dif W_s
$$
and
\begin{align*}
\eta_t=M_t\left[\zeta_0+\int^t_0 M^{-1}_s(\zeta^{(1)}_s-\<\bar\alpha_s,\zeta_s^{(2)}\>)\dif s+\int^t_0M^{-1}_s \zeta^{(2)}_s\dif W_s\right].
\end{align*}
Hence, 
\begin{align*}
\|\eta^*_T\|_{q_0}&\leq\|M^*_T\zeta_0\|_{q_0}+\left\|M^*_T(M^{-1})^*_T\int^T_0|\zeta^{(1)}_s|\dif s\right\|_{q_0}\\
&\quad+\left\|M^*_T(M^{-1})^*_T\int^T_0|\bar\alpha_s|\cdot|\zeta_s^{(2)}|\dif s\right\|_{q_0}\\
&\quad+\left\|M^*_T\sup_{t\in[0,T]}\left|\int^t_0M^{-1}_s\zeta_s^{(2)}\dif W_s\right|\right\|_{q_0}\\
&=:I_1+I_2+I_3+I_4.
\end{align*}
Noticing that by \eqref{C2},
\begin{align}\label{C3}
|\bar\beta_s|\leq|\beta_s|,\ \ |\bar\alpha_s|\leq|\alpha_s|,
\end{align}
for any $p\in\mR$, by \eqref{C1}, H\"older's inequality and Doob's maximal inequality, we have
\begin{align}\label{Expo}
\mE\left(\sup_{t\in[0,T]}M_t^p\right)<\infty.
\end{align}
Thus, by H\"older's inequality and \eqref{Expo}, we have
$$
I_1\leq C\|\zeta_0\|_{q_1},\ \ 
I_2\leq C\left\|\int^T_0|\zeta^{(1)}_s|\dif s\right\|_{q_2},
$$
and by \eqref{C3},
\begin{align*}
I_3&\leq\left\|M^*_T(M^{-1})^*_T\left(\int^T_0|\alpha_s|^2\dif s\right)^{1/2}\left(\int^T_0|\zeta_s^{(2)}|^2\dif s\right)^{1/2}\right\|_{q_0}\\
&\leq C\left\|\left(\int^T_0|\zeta_s^{(2)}|^2\dif s\right)^{1/2}\right\|_{q_3}=C\left\|\int^T_0|\zeta_s^{(2)}|^2\dif s\right\|^{1/2}_{q_3/2}.
\end{align*}
Similarly, by H\"older and Burkholder's inequalities, we also have
$$
I_4\leq C\left\|\int^T_0|\zeta_s^{(2)}|^2\dif s\right\|^{1/2}_{q_3/2}.
$$
Combining the above estimates, we obtain \eqref{HG8}.
\end{proof}

Let $f$ be a locally integrable function on $\mR^d$.
The Hardy-Littlewood maximal function is defined by
$$
\cM f(x):=\sup_{0<r<\infty}\frac{1}{|B_r|}\int_{B_r}f(x+y)\dif y,
$$
where $B_r:=\{x\in\mR^{d}: |x|<r\}$.
The following result can be found in  \cite[Appendix A]{Cr-De}.
\bl\label{Le2}
(i) There exists a constant $C_d>0$ such that for all $f\in C^\infty(\mR^{d})$ and $x,y\in \mR^{d}$,
\begin{align}
|f(x)-f(y)|\leq C_d |x-y|(\cM|\nabla f|(x)+\cM|\nabla f|(y)).\label{Es2}
\end{align}
(ii) For any $p>1$, there exists a constant $C_{d,p}$ such that for all $f\in L^p(\mR^d)$,
\begin{align}
\|\cM f\|_p\leq C_{d,p}\|f\|_p.\label{Es30}
\end{align}
\el

\end{document}